\documentclass[
useAMS,
usenatbib,
usegraphicx]
{tellus}

\usepackage{color}
\usepackage{rotating}
\usepackage{multirow}
\usepackage{amsmath,amssymb}
\usepackage{subfigure} 
\usepackage{hyperref}

\newcommand{\xvec}{\ensuremath{\mathbf{x}}}
\newcommand{\yvec}{\ensuremath{\mathbf{y}}}

\begin{document}

\title{Predicting flow reversals
  in chaotic natural convection using
  data assimilation}

\author[]{By Kameron Decker Harris$^1$\thanks{Corresponding
    author.\hfil\break e-mail: kameron.harris@uvm.edu}, 
  El Hassan Ridouane$^1$,
  Darren L. Hitt$^2$, 
  and Christopher M. Danforth$^1$} 
\affiliation{
  $^{1}$Department of Mathematics and Statistics,
  Complex Systems Center, 
  \& the Vermont Advanced Computing Center,\\
  University of Vermont, Burlington, Vermont 05405, USA\\
  $^{2}$Department of Mechanical Engineering,
  University of Vermont, Burlington, Vermont 05405, USA
} 

\history{\today}
\maketitle

\begin{abstract} 
  A simplified model of natural convection, 
  similar to the Lorenz (1963) system, is compared to
  computational fluid dynamics simulations of a thermosyphon
  in order to test 
  data assimilation methods and better understand
  the dynamics of convection.
  The thermosyphon is represented by a long time flow simulation,
  which serves as a reference ``truth''. 
  Forecasts are then made using the Lorenz-like
  model and synchronized to noisy and limited observations of 
  the truth using data assimilation.
  The resulting analysis is observed to infer dynamics absent from the
  model when using short assimilation windows.

  Furthermore, chaotic flow reversal occurrence and residency times in
  each rotational state are forecast using analysis data.
  Flow reversals have been successfully forecast in the 
  related Lorenz system, as part of a perfect model experiment,
  but never in the presence of significant model error or 
  unobserved variables.
  Finally, we provide new details concerning the fluid dynamical processes
  present in the thermosyphon during these flow reversals.
\end{abstract}


\section{Introduction}

Forecasting methodologies, traditionally motivated by numerical
weather prediction (NWP), can find applications in other
fields such as engineering \citep{savely1972a}, 
finance \citep{sornette2006a, bollen2011a}, 
epidemiology \citep{ginsberg2009a}, and marketing \citep{asur2010a}.
Techniques borrowed from the weather forecasting community may
prove to be powerful for forecasting these other types of complex systems.
Fluid systems can be particularly challenging
due to dynamics taking place at multiple interacting spatial and
temporal scales. However, because of their relationship to NWP, 
fluid systems are among the most studied in the context of forecasting.

In this paper, we show that the flow in a 
computational fluid dynamics (CFD) simulated thermosyphon
undergoing chaotic convection can be accurately forecast using
an ordinary differential equation (ODE) model akin to the classic
\cite{lorenz1963} system.
The thermosyphon, a type of  
natural convection loop or non-mechanical heat pump,
can be likened to a toy model of climate.
Thermosyphons are used in solar water
heaters \citep{belessiotis2002}, cooling systems for computers
\citep{beitelmal2002a}, roads and railways that cross permafrost
\citep{lustgarten2006}, nuclear power plants \citep{detman1968a, beine1992,
  kwant1992}, and other industrial applications. 
In these heat pumps, buoyant
forces move fluid through a closed loop,
and at high amounts of forcing they can exhibit complex aperiodic behavior.
As first suggested by \cite{lorenz1963}, this
is illustrative of the unpredictable convection
behavior observed in weather and climate dynamics.

Synthetic observations of the thermosyphon are combined with model data
to produce new forecasts in the process known as data assimilation (DA).
DA is a generic method of combining observations
with past forecasts to produce the analysis, an approximately
optimal initial condition (IC) for the next forecast cycle. 
Another interpretation of the analysis is that it is a ``best guess''
for the true system state as represented in the phase space of the model.
DA can be used as a platform for the reanalysis of past observations, 
in which the dynamical model plays a key role in constraining
the state estimates to be physically realistic \citep{compo2006a}.

In the present study, we use an analysis of simulated 
thermosyphon mass flow rate data to
explain the heat transport processes occurring
during chaotic flow reversals, and to 
inform empirical forecasts of the occurrence of these flow reversals.
Although the flow reversals are chaotic, we show they have
short-term predictability, quantifying the extent to which this
is possible with our methods.

This paper is structured in the following way: 
In Sec.~\ref{sec:model}, 
we explain the CFD simulation used to generate a
synthetic true state or ``nature run'' of the thermosyphon
and the separate forecasting model.
In Sec.~\ref{sec:da} we present an overview of how DA was 
applied to this experiment and its performance.
In Sec.~\ref{sec:reversals} we explain and present the results for
flow reversal and rotational state residency time forecasts.
Finally, Sec.~\ref{sec:conclusion} contains concluding remarks.
In Appendices S1-4 in the Supporting Information we present a
derivation of the model,
detail the tuning of model parameters, 
and explain in detail the DA methods used.

\begin{figure} \centering
  \includegraphics[width=.9\linewidth]{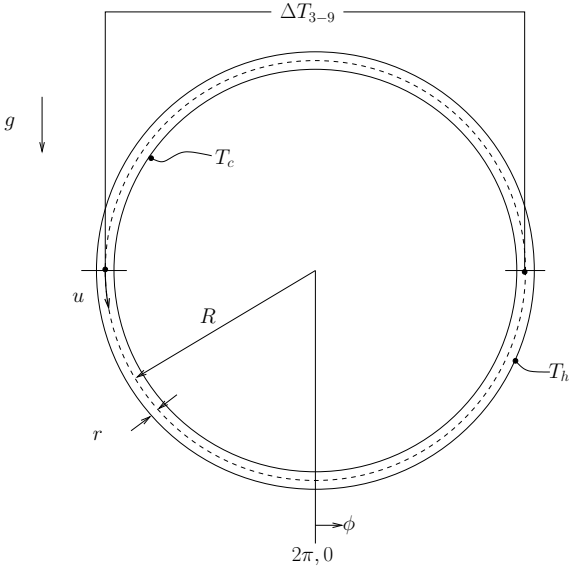}
  \caption{The thermosyphon has a simple circular geometry. The bottom
    wall is heated to a constant hot temperature $T_h$ while the top wall
    is maintained at the temperature $T_c$, creating a
    temperature inversion of hot fluid below cold fluid. If conduction alone
    cannot stabilize this temperature inversion, then the fluid will begin
    to rotate and convection becomes the dominant process of heat
    transfer. The relevant model state variables are proportional to
    the bulk fluid velocity $u$ and the temperature difference across the loop
    $\Delta T_{3-9}$. For counterclockwise flow, as indicated by the arrow
    near 9 o'clock, fluid velocity $u > 0$ and
    $\Delta T_{3-9}$ is typically $> 0$.
    The radius ratio $R/r=24$ used in our experiments is
    shown to scale.}
  \label{fig:loopcartoon}
\end{figure}

\section{Models and the data assimilation algorithm}
\label{sec:model}

Following previous experiments that examined the 
periodic \citep{keller1966} and chaotic 
\citep{welander1967,
  creveling1975, gorman1984, gorman1986, ehrhard1990, yuen1999a,
  jiang2003,burroughs2005, desrayaud2006,yang2006a, ridouane2009}
behavior of toroidal thermosyphons, we also consider a circular
thermosyphon geometry. Picture a vertically-oriented hula
hoop, as shown in Fig.~\ref{fig:loopcartoon}. 
An imposed wall temperature $T_h$ on the lower half of the loop
($-\frac{\pi}{2} < \phi < \frac{\pi}{2}$) heats the fluid contained
in this section. Similarly, a wall temperature $T_c < T_h$ is
imposed on the upper half
($\frac{\pi}{2} < \phi < \frac{3\pi}{2}$) to cool the upper section
(Fig.~\ref{fig:loopcartoon}).
The forcing, proportional to the temperature difference $\Delta T_w =
T_h-T_c$, is constant. We focus on the case of developed flow,
ignoring transient behavior.

The behavior of the fluid can be qualitatively understood as follows.
As the heating parameter is increased, the 
flow behavior transitions from a conduction state 
(conducting equilibrium) 
to a steady, unidirectional state of convection 
(convecting equilibrium). 
No particular rotational state 
(clockwise, CW, or counterclockwise, CCW)
is favored due to symmetry. At still higher heating values,
chaotic flow oscillations can be observed.
In the chaotic regime, the flow is observed to oscillate around
one unstable convecting equilibrium state until flow reversal.
Each flow reversals causes the system to transition
between CW and CCW rotational states.

\subsection{Thermosyphon simulation} 
\label{sec:thermosyphon}

The reference state of the thermosyphon is
represented by a CFD-based numerical
simulation in two spatial dimensions (2D). 
The details of the computational 
model have been described in detail in a previous study by
\cite{ridouane2009}; however, for completeness, 
we summarize here its essential elements. 

It is assumed that the temperature differential $\Delta T_w$
is sufficiently small so that temperature-dependent
variations of material properties can be regarded
as negligible, save for the density.  
The standard Boussinesq 
approximation is invoked and all fluid properties are assumed to be 
constant and evaluated at the reference temperature $(T_h+T_c)/2$.  
The flow is assumed to be laminar, two-dimensional,
with negligible viscous dissipation due to low velocities. 
Under these circumstances, the governing
dimensionless equations are the unsteady, 2D laminar Navier-Stokes 
equations along with the energy equation and equation of 
state for the density. 
No slip velocity boundary conditions are
imposed on the walls and isothermal boundary conditions of $T_h$ and
$T_c$ are imposed on the heated and cooled 
lower and upper walls, respectively.

The dimensionless control parameter for convection is the
Rayleigh number, defined here as 
\begin{equation}
  \mathrm{Ra} = \frac{ 8 g \gamma r^3  \Delta T_w}{\nu \kappa} ,
\end{equation}
where $g$ is the gravitational acceleration, $\gamma$ is the thermal 
expansion coefficient, $\nu$ is the kinematic viscosity, and $\kappa$ is 
the thermal diffusivity.

The one dimensionless geometric parameter is the ratio
of major (loop) radius $R$ to minor (tube) radius $r$, hereafter
referred to as the radius ratio.
Consistent with the previous study,
the dimensions 
of the loop are chosen with $R$ = 36 cm and $r$ = 1.5 cm to yield a 
radius ratio of 24. 

As in the classic Rayleigh-B\'{e}nard problem, the Rayleigh number 
determines the onset of convection in the thermosyphon.  For the 
numerical simulations on this fixed geometry, a range of Rayleigh 
numbers can be imposed by varying the value of the gravitational 
acceleration. As the Rayleigh number is increased from zero, the 
flow behavior transitions  from a stationary, conduction state to a 
steady, unidirectional state of convection. At still higher values of 
Ra, chaotic flow oscillations can be observed.  Unless otherwise indicated, 
the simulation results presented in this paper correspond to a value of 
Ra= $1.5 \times 10^5$, which is within the chaotic regime.

All numerical simulations were performed using the commercial
CFD software \cite{fluent}, which is based on
the finite-volume method. 
(An example of the output is shown in
Fig.~\ref{fig:regchangeprofiles}, in the discussion of flow reversals.)
During the course 
of the simulations, the
time-varying mass flow rate, a scalar denoted by
$q$ and proportional to $u$,
is saved at 10 s intervals. This reporting interval is conservative, as
laboratory thermocouples can be sampled more than once pre second. 
In doing so, a 
time series of the ``true'' synthetic thermosyphon state is recorded
to be used in a forecasting scheme.

\subsection{Forecast model}
\label{sec:EM}

The Ehrhard-M\"uller (EM) system is a three-variable ODE derived
specifically to model bulk flow in
the thermosyphon (\citealp{ehrhard1990}; also see Appendix S1 in the
Supporting Information for an alternative derivation).
Written in dimensionless
form, the governing equations are
\begin{subequations} \label{eq:EM}
  \begin{align} 
    \frac{dx_1}{dt'} &= \alpha \left(x_2 - x_1\right)\\
    \frac{dx_2}{dt'} &= \beta x_1 - x_2 \left( 1+K h(|x_1|) \right) - 
    x_1 x_3\\ 
    \frac{dx_3}{dt'} &= x_1 x_2 - x_3 \left( 1+K h(|x_1|) \right) .
  \end{align}
\end{subequations} 
The state variable $x_1$ is proportional to the mass flow rate or
mean fluid velocity, $x_2$ to
the temperature difference across the convection cell ($\Delta
T_{3-9}$, measured between 3 o'clock and 9 o'clock), and $x_3$ 
to the deviation of the vertical temperature profile
from the value it takes during conduction; specifically,
$x_3 \propto (\frac{4}{\pi} \Delta T_w - \Delta T_{6-12})$,
where $\Delta T_{6-12}$ is the temperature difference
measured between 6 o'clock and 12 o'clock.
The parameter $\alpha$ is comparable to the Prandtl number, 
the ratio of momentum diffusivity
and thermal diffusivity. 
Similar to the Rayleigh number, the heating
parameter $\beta \propto \Delta T_w$ 
determines the onset of convection as well as the
transition to the chaotic regime.  
Finally, $K$ determines the magnitude of variation of the
wall heat transfer coefficient with velocity. 
The functional
form of that variation is determined by 
$h: \mathbb{R}^{+} \to \mathbb{R}^{+}$, where
\begin{equation}
  \label{eq:H}
  h(x) = \left\{
    \begin{array}{lr}
      \frac{44}{9} x^2 - \frac{55}{9} x^3 + \frac{20}{9} x^4 & 
      \text{when } x < 1 \\
      x^{1/3} & \text{when } x \geq 1
    \end{array}
  \right. .
\end{equation}
The interested reader is referred to Appendix S1 in the Supporting 
Information for an explanation of this piecewise form, which differs
slightly from the original model of \cite{ehrhard1990}.

Note that when $K=0$, the system is analogous to the \cite{lorenz1963}
system with geometric factor (Lorenz's $b$) equal to one.
The lack of a geometric factor in the EM system is due to
the circular geometry of the convection cell. 
Lorenz equations have been widely used in nonlinear dynamics to study
chaos and in NWP as a model system for testing DA
\citep{miller1994a,yuen1999a,annan2004,evans2004,yang2006a,kalnay2007a}.

When in the chaotic parameter regime, the EM system
exhibits growing oscillations in the $x_1$ and $x_2$
state variables around their convecting equilibrium values until 
flow reversal. In this system, the CCW rotational
state is characterized by $x_1 > 0$ and $x_2 > 0$, and the CW
rotational state by $x_1 < 0$ and $x_2 < 0$.
However, one should note that
near a flow reversal $x_1$ and $x_2$ can 
have opposite signs, because zero-crossings of the
$x_1$ variable typically lag behind those of $x_2$.

The parameters found to match the simulated thermosyphon
were $\alpha = 7.99$, 
$\beta = 27.3$, and $K = 0.148$.
The characteristic time and mass flow rate scales, 
used to transform the dimensionless model variables
$t'$ and $x_1$ into dimensional time and ``observations'' 
of mass flow rate, were 631.6 s
and 0.0136 kg/s,
respectively. 
The $q$ scale is the one nonzero
entry in the observation operator $\mathbf{H}$, Eqn.~\eqref{eq:obsop}. 
The above parameters were found using a multiple shooting
algorithm explained in Appendix S1.2 in the Supporting Information.
Numerical integration of this
autonomous ODE was performed with a fourth-order Runge-Kutta
method and timestep 0.01 (corresponding to 6.316 s)
in \cite{matlab}.

\subsection{Data assimilation}
\label{sec:dadesc}

\begin{figure}
  \centering 
  \includegraphics[totalheight=0.24\textheight,
  width=\linewidth,
  keepaspectratio=true]{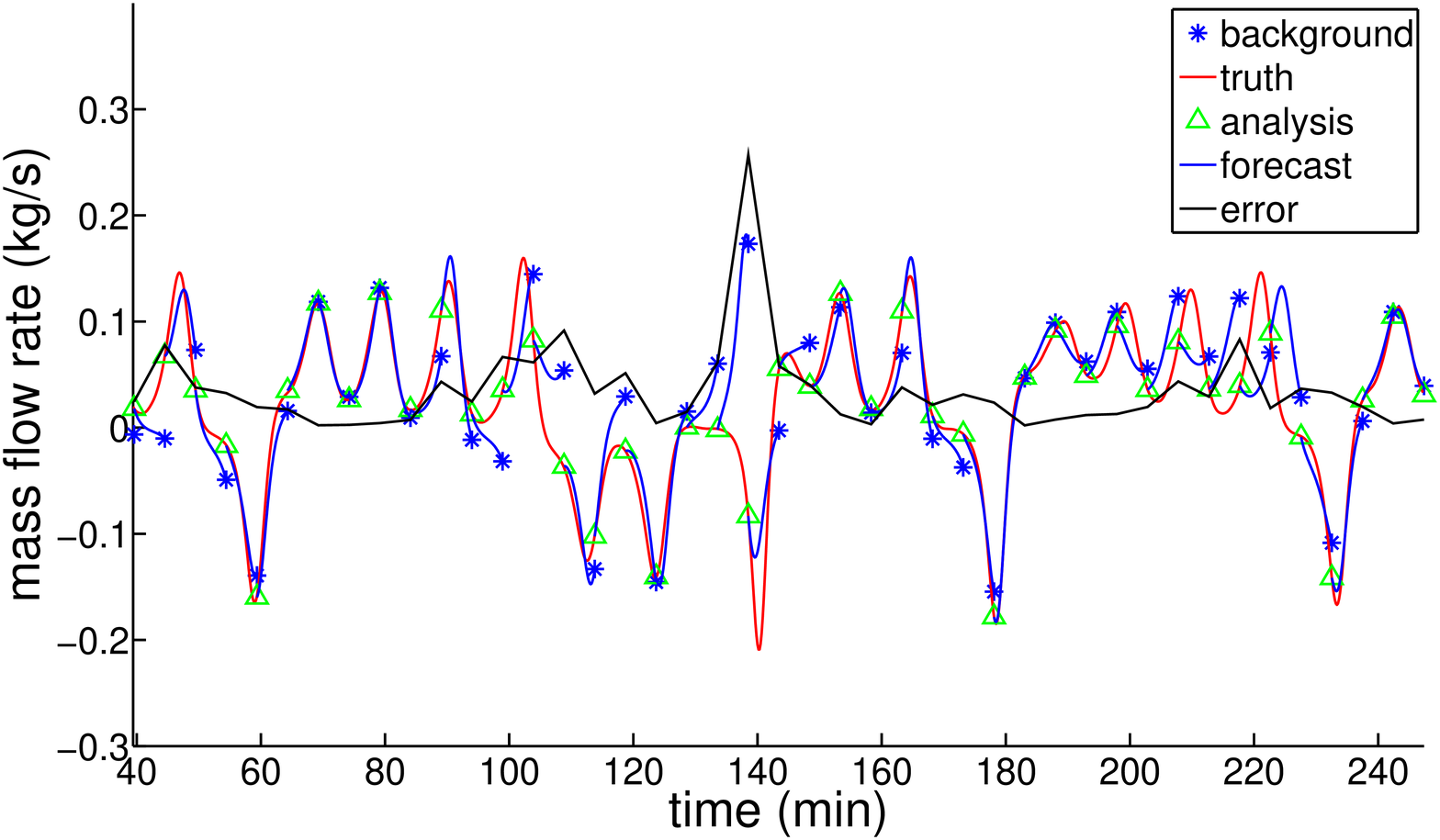}
  \caption{An illustration of the basic
    predict-observe-update DA cycle. The EM model states
    (background forecast and analysis)
    are transformed into observations of the mass flow rate $q$
    by the observation operator (Eqn.~\eqref{eq:obsop})
    for comparison with the truth.
    Here, using the 3D-Var algorithm and
    an assimilation window of 5 min, the
    filter has satisfactory overall performance 
    (scaled error $\approx$ 35\%).
    Note the error spike around 135 min
    when the forecast and truth end up in 
    different rotational states. 
    The largest errors tend to occur at or near flow reversals due
    to inherent sensitivity near that transition and to
    the qualitatively different behavior of the different flow directions.
  }
  \label{fig:daexample}
\end{figure}

DA is the process by which observations of a dynamical system are
combined with forecasts from a model to estimate error covariances and
calculate an optimal estimate for the current state of the system,
called the analysis. The inherent difficulties
are compounded by the fact that the forecaster uses an inexact
forecasting model and never knows the true state of the dynamical
system.
The number of state variables in a NWP model is typically
$\mathcal{O}(10^3)$ times larger than the 
number of observations. Nevertheless,
the analysis becomes the IC for a new forecast.
The time interval between successive applications of the DA algorithm, 
i.e. the time between analysis steps
(usually determined by the availability of observations but here
allowed to vary), 
is called the assimilation window.
The process is illustrated in Fig.~\ref{fig:daexample}.

A variety of filters are capable of solving the DA problem. 
The canonical example is the Kalman filter (KF; \citealp{kalman1960}),
the optimal state estimation algorithm for a linear system. 
One of DA's first applications was to trajectory estimation and correction of
missiles and rockets \citep{savely1972a}. 
A number of nonlinear DA schemes are implemented in this study.
In 3D variational DA 
(3D-Var; here 3D refers to the spatial dimensions for weather models), 
the background error covariance is estimated a single time, offline, prior
to the data assimilation procedure. In the extended Kalman filter
(EKF), background error is evolved according to the linear tangent model,
which approximates the evolution of small perturbations about the
trajectory. 
Ensemble Kalman filters (collectively EnKFs) 
use ensembles of forecasts to estimate the
background error and better capture nonlinear behavior. 
The methods examined in this study were 3D-Var, the EKF,
the ensemble square root Kalman filter (EnSRF), and the 
ensemble transform Kalman filter (ETKF).
Detailed descriptions of each method
are included in Appendix S2 in the Supporting Information.
A full review of DA is beyond the scope of the present paper; 
for a comprehensive treatment, we refer the reader to \cite{kalnay2002}.

\section{Data assimilation experiments}
\label{sec:da}

\subsection{Methods}
\label{sec:da.setup}

A perfect model experiment, in which the Lorenz equations were used to 
forecast a synthetic truth created by the exact same system,
was tested first but not included here. 
We found analysis errors similar to
those reported by \cite{yang2006a} (3D-Var and EKF) and
\cite{kalnay2007a} (ETKF), using the same model and tuning
parameters. This ensured that the DA algorithms were working
before applying them to the synthetic thermosyphon data.

As stated in Sec.~\ref{sec:thermosyphon}, forecasts
of the thermosyphon are made
observing one scalar variable, the mass flow rate 
$q \propto x_1$.
Gaussian noise with standard
deviation equal to $6\times 10^{-4}$ kg/s, approximately $0.8\%$ of
the mass flow rate climatological mean, 
$\sqrt{\langle q^2 \rangle} = 0.075812$ kg/s, is
added to the synthetic truth to create observations.
The relative magnitude of this error is
comparable to that of experimental measurements.

The EM model is used in the forecast step to
integrate the analysis forward in 
time and create the new background forecast.
The end results of applying DA are a background and analysis
timeseries of $x_1, x_2, x_3$, informed by both
the timeseries of thermosyphon mass flow rate and 
the EM model dynamics.

In this realistic forecasting scenario, where only limited 
information about the true state is available,
the observations of state variable $q$
provide the only validation. 
For this reason, we calculate the forecast errors 
in observation space.
These are given as root mean square error (RMSE),
where $\mathrm{RMSE} = \sqrt{\langle \delta q^2 \rangle}$. 
The residual at a
specific assimilation cycle is given by $\delta q = q -
\mathbf{H}\xvec^b$. Here, $\xvec^b$ is the background forecast made
by the model,
and $\mathbf{H}: \mathbb{R}^3 \to \mathbb{R}$ is the
linear observation operator
\begin{equation}
  \label{eq:obsop}
  \mathbf{H} = \left[ 0.0136, 0, 0 \right]
\end{equation} in units of kg/s.
All errors are then scaled by $\sqrt{ \langle q^2\rangle }$,
the climatology of $q$.
Analysis error is a common metric for assessing DA performance
in perfect model experiments. In this study, however, we
assert that background error is preferable.
Analysis error in observation space,
which will be small even for large assimilation windows, 
is not an appropriate
metric for assessing model performance since it can disagree
substantially with the background error. 
For example, 3D-Var in one experiment
with a 10 minute assimilation window yielded analysis and background
scaled errors of 0.08 and 0.86, respectively. 
The analysis error would seem to indicate that forecasting is
doing a good job, but the background error shows that background
forecasts are essentially meaningless. 
The filter, however, accounts for this and weights the
observations heavily over the background
forecasts when producing the analysis. 
Since we are concerned with forecasting, 
background error is a more representative metric.

When applying DA to nonlinear systems, some type of covariance inflation
is performed to prevent filter divergence due to error underestimation. 
\cite{kalnay2007a}
found that a Lorenz forecasting model with a slightly different forcing
parameter required a 10-fold increase in the
multiplicative inflation factor when using a 3 member EnKF.
Model error is more pronounced for our forecasts,
since the EM model is a reduced approximation of the numerically simulated
thermosyphon.
We relied upon additive and multiplicative background
covariance inflation to capture model error. 
Additive inflation was particularly important for the stability
of the EKF and EnKFs. 
Additive noise provides
a different exploration of dynamically accessible regions of state
space, and it would be interesting to explore why additive versus 
multiplicative is preferred in certain cases, although this is beyond
the scope of this paper.
The specifics of how inflation was performed 
and tuned, and the parameters used are given in
Appendix S2 in the Supporting Information.

All EM and DA parameter tuning was performed using a separate
mass flow rate time series than was used for validation. 
Each DA algorithm was allowed 500 cycles to spin-up, and its
performance was measured over the following 2500 cycles.
Ensemble size in each case was set to 10 members.

\subsection{Results}
\label{sec:da.results}

\begin{figure}
  \centering
  \includegraphics[width=\linewidth]{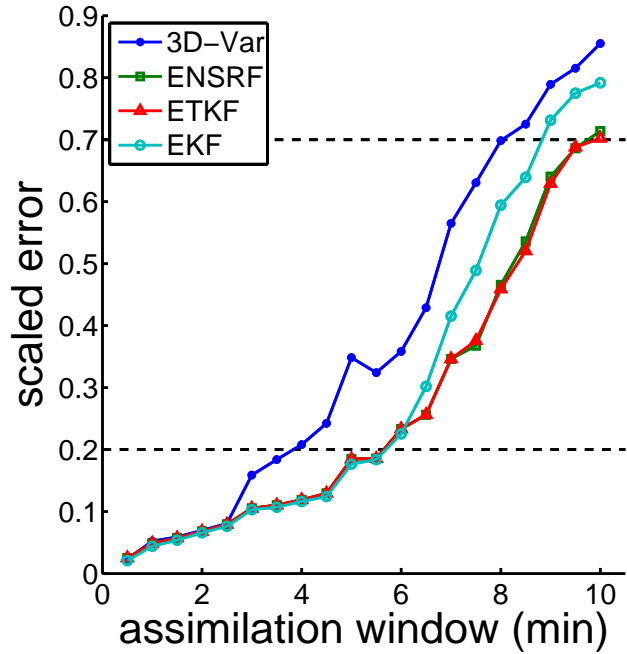}
  \caption{Background RMSE scaled by 
    $\sqrt{\langle q^2 \rangle}$ 
    over 2500 assimilation cycles
    plotted for different DA
    algorithms and varying assimilation windows. 
    As the window becomes larger,
    the error increases towards saturation. 
    The lower dashed line
    in the main figure shows the limit of a ``perfect'' forecast while the
    upper demarcates a ``useless'' forecast.
  }
  \label{fig:EMFwindows}
\end{figure}

With proper tuning, all DA algorithms were capable of synchronizing
the EM model to observations of mass flow rate alone. 
As the assimilation window
increased, scaled background error increased in a sigmoidal fashion,
as expected (see Fig.~\ref{fig:EMFwindows}).
For assimilation windows up to 2.5 min, all DA algorithms
have nearly indistinguishable errors. 
For assimilation windows between
3 and 6 min, 3D-Var performs noticeably worse than the other
methods which remain indistinguishable. Then, with assimilation windows
greater than 6 min, the ensemble methods (EnSRF and ETKF) outperform 
the EKF noticeably.
This is perhaps surprising, at first glance,
because the ensemble
size is significantly smaller than the dimension
of the simulated thermosyphon state space ($\mathcal{O} (10^5)$ variables).
However, we know the thermosyphon dynamics effectively take place
on the EM equations' attractor (a manifold in three dimensions).
The superior performance of EnKFs here is likely due
to the ensemble methods capturing nonlinear
effects which dominate at larger windows.

Following the historical S1 score convention, scaled error
above 70\% is considered a ``useless'' forecast, while under 20\% the
forecast is ``perfect'' 
\citep{kalnay2002}.
Perfect forecasts for 3D-Var were found
up to a 4 minute assimilation window, 
while the other methods (EnSRF, ETKF, and EKF) produced perfect
forecasts with assimilation windows 1.5 minutes longer.

A persistent spike in background error for the 5 minute assimilation window
(Fig.~\ref{fig:EMFwindows})
is possibly due to that time
being approximately the same as the characteristic
period of oscillations in $q$ 
(evident in Fig.~\ref{fig:daexample}). 
We conjecture that this
may lead to a type of resonance in the DA-coupled EM system
which degrades DA performance.

\begin{figure} 
  \centering 
  
  \includegraphics[totalheight=.65\textheight, 
  width=\linewidth, keepaspectratio=true]{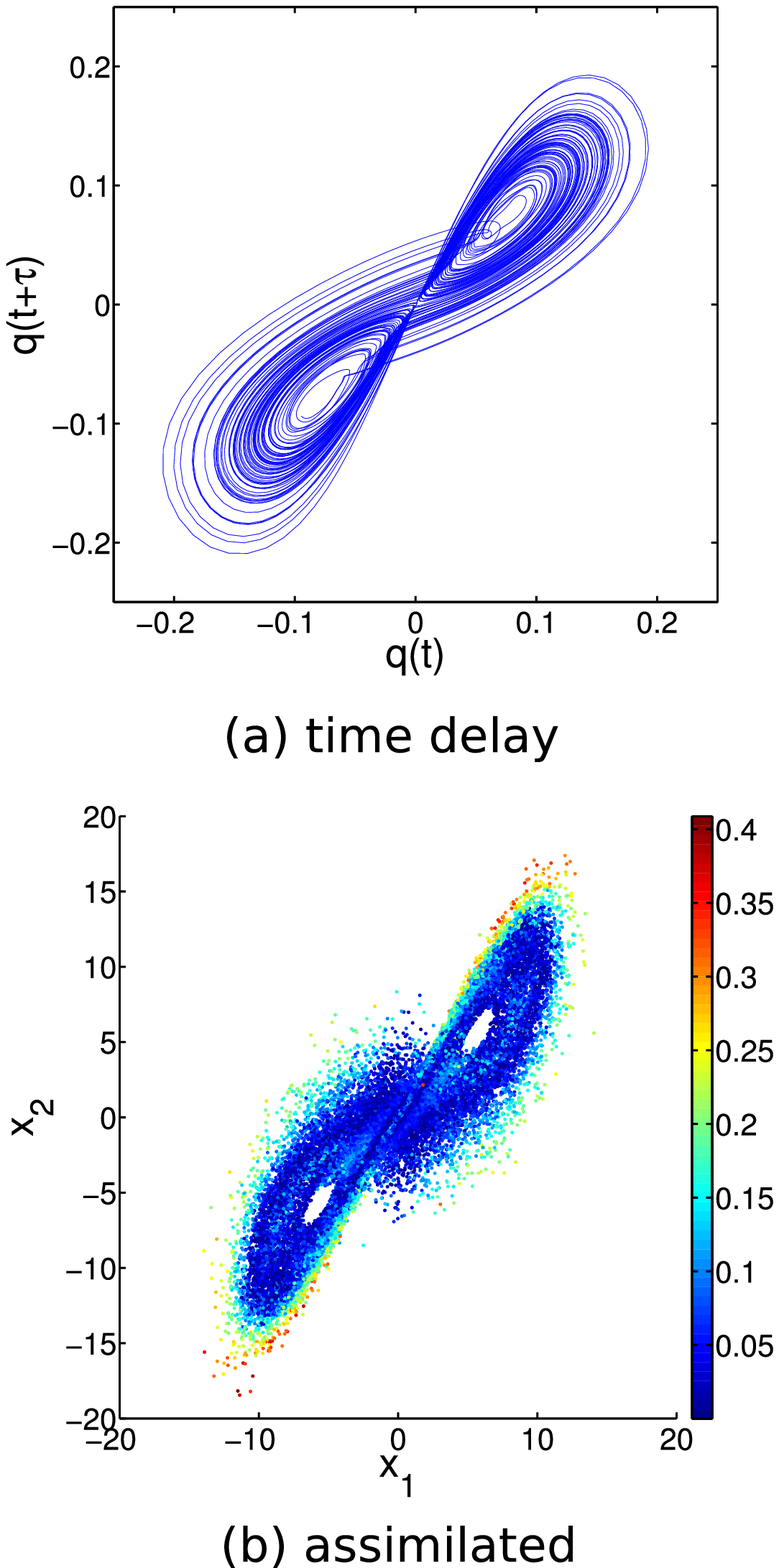}

  \caption{Two views of the numerically simulated thermosyphon attractor. 
    A time-delay reconstruction, using the monitored mass flow rate, is
    shown in (a). In (b),
    plotted points show $x_1$ and $x_2$ of the EM analysis generated by EKF
    with an assimilation window of 120 s. Each is colored by the
    scaled background forecast error at that point. 
    The delay time of 60 s used for (a) was chosen because it yielded a good
    approximation of the attractor in (b). Note how in (a) trajectories that
    move through the far edge of either lobe create distinctive loops near
    the center of the opposite lobe. This is an example of dynamics
    which are not present in the EM model without DA.
    It may explain
    the higher error for points in (b) at the far
    edge of each attractor lobe.
    See text for further description and Fig.~\ref{fig:perfectstorm}
    for another example.} 
  \label{fig:attractor}
\end{figure}

Besides these results pertaining to forecast skill,
we also found that the DA algorithms infer thermosyphon dynamics 
which are absent from the EM model. 
In Fig.~\ref{fig:attractor} we
see the simulated thermosyphon's attractor obtained by both a
time-delay embedding 
(Fig.~\ref{fig:attractor}(a); \citealp{alligoodchaosbook}) 
and a projection of the EM analysis states
to the $x_1$-$x_2$ plane (Fig.~\ref{fig:attractor}(b)).
If the thermosyphon fluid flow stalls in the midst of a reversal, 
fluid in the bottom can quickly heat up while
that in the top is cooled,
leading to an unstable, strong temperature inversion.
This causes the fluid to move very quickly in the reversed direction, 
but this new direction also ends up being unstable,
and a new flow reversal can occur immediately.
Absent DA, the EM model system does not exhibit this behavior 
of stalling followed by large swings of the trajectory.

In the time-delay embedding (Fig.~\ref{fig:attractor}(a)), 
this phenomenon is exhibited by small loops in the trajectory as
it moves near the convective fixed points. 
The flow stalls when the system state
is near the conductive fixed point at the origin,
then it swings wildly which brings it near the convective fixed
point, but in such a way that it does not end up spiraling
outward in the usual fashion as during a normal
flow reversal, as exhibited by the Lorenz equations. 
Instead, it quickly reverses again, which we call non-Lorenz behavior.
This non-Lorenz behavior is further elaborated upon in 
Sec.~\ref{sec:revoccur.new}.
Forecast skill is worst at the far edges of the
assimilated attractor (Fig.~\ref{fig:attractor}(b)). 
This could be due to the wild swings of the EM trajectory
after being ejected from the region of state space
near the conducting equilibrium, or to the nonlinear dynamical
instabilities at the edge of the attractor found by
\cite{palmer1993a} and \cite{evans2004}.

\begin{figure}
  \centering
  \includegraphics[width=\linewidth]{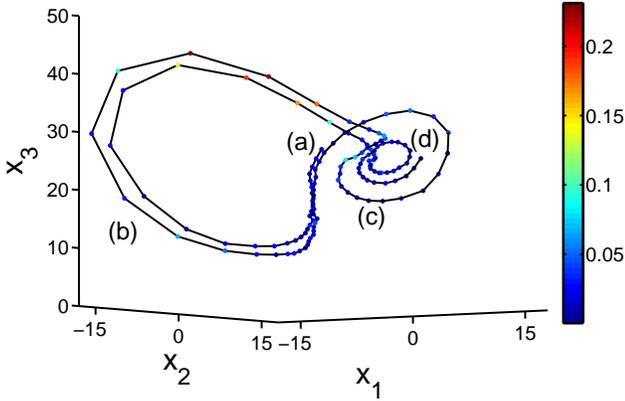}
  \caption{The assimilated trajectory during a very
    large-oscillation, non-Lorenz flow reversal.
    The EKF algorithm with a 30 s assimilation window was used.
    Color indicates the scaled background error.
    The state starts (a) in
    the CCW rotational region and stalls near the conducting
    state for an extended period, causing fluid in the bottom to heat up,
    manifesting in an $x_3$ that creeps towards 0 with 
    $x_1$,$x_2 \approx 0$, before the state swings 
    quickly (b) through one oscillation in the CW
    rotational state. This is followed by one oscillation in the  
    CCW state (c) before another stall near the conducting point and
    subsequent swing (b again) before settling into Lorenz-like
    oscillations (d). Note that the filter only observes the $x_1$
    variable but is able to reconstruct the dynamics in the
    full state space.
  }
  \label{fig:perfectstorm}
\end{figure}

We also explicitly show one of these stalled flow reversals
in Fig.~\ref{fig:perfectstorm}, where we plot the EKF-assimilated
EM trajectory using a 30 s assimilation window.
When the fluid stalls,
the $x_3$ variable moves closer to 0 (i.e. $\Delta T_{6-12}$ increases)
while $x_1$ and $x_2$ (proportional to $q$ and $\Delta T_{3-9}$, respectively) 
are approximately 0, reflecting the growing temperature inversion 
while the fluid remains stationary. When the fluid starts to move,
the assimilated trajectory swings wildly to the left attractor lobe
(CW rotational state), 
then right (CCW rotational state). 
The trajectory undergoes another stall-swing cycle
before finally resuming Lorenz
behavior, where the trajectory spirals outward from the
CCW convecting equilibria.
This contrasts the Lorenz and EM model dynamics, 
for which large deviations
in the system state from convecting equilibrium are driven
close to the other convecting equilibrium during a flow reversal,
which stabilizes the system. 
See also Sec.~\ref{sec:revoccur.new},
 Fig.~\ref{fig:EMFregchange}, and the accompanying discussion.
This result remains unchanged for the other DA algorithms also using a 30 s 
assimilation window. The inference remains
using EKF and a 60 s assimilation window, 
but the trajectory appears much noisier, leading us to
believe that this is due to the rapid update. 
With larger assimilation windows, the trajectory becomes
uninterpretable as error in the unobserved variables increases.

\section{Flow reversal experiments}
\label{sec:reversals}

\subsection{Experimental setup}
For the purpose of flow reversal forecasts,
we picked a single DA algorithm and assimilation window.
In this Section, all analyses were generated by the extended Kalman filter
and an assimilation window of 30 s. This interval corresponds to 
5 time steps of the model and is shorter than that used in 
\cite{yang2006a} and \cite{kalnay2007a}.
The following could certainly
be repeated using other algorithms, observations, and assimilation windows,
but this was beyond the scope of this paper.
The flow reversal tests in Sec.~\ref{sec:revoccur.tests}
and the residency time forecasts in 
Sec.~\ref{sec:revoccur.dur}
were tuned and validated on separate analysis timeseries.
The length of the tuning and validation timeseries were approximately
39 and 93 days, respectively.

\subsection{Occurance of flow reversals: traditional explanation}
\label{sec:revoccur.trad}
The first explanation of the mechanism responsible for
flow reversals was presented by \cite{welander1967} and repeated by
\cite{creveling1975}. Welander, who was also the first to discover
that thermosyphons exhibit aperiodic oscillatory behavior, explained
the instability of steady convecting flow by considering a thermal
anomaly or ``warm pocket'' of fluid. For low heating rates, the
convecting equilibrium is stable because viscous and thermal
dissipation are in phase, thus an increase (decrease) in flow rate
leads to an increase (decrease) in friction and a decrease (increase)
in buoyancy, and such perturbations are damped out. At higher heating
rates, the warm pocket is amplified with each cycle through the loop
due to out of phase viscous and thermal dissipation's. Welander
explained that when the warm pocket emerges from the heating section
and enters the cooling section, it feels a greater buoyant force than
the surrounding fluid and accelerates, exiting the cooling section
quickly, giving it less time to radiate away its energy. As the pocket
moves into the section with warm boundary, 
the buoyant force it experiences is again
higher than normal, so now the pocket decelerates and passes slowly
through the heating section, gaining more energy.
This positive feedback effect causes the pocket to grow hotter and
larger with each pass through the loop.
These oscillations in the fluid temperature and velocity 
do not grow unhindered, however.
The pocket eventually becomes large and hot enough that its
descent towards the heating section is stopped entirely by its own
buoyancy. Without movement, the pocket dissipates, but its remnant
heat biases new rotation in the opposite direction, and the flow
reverses.

In the Lorenz and EM systems, 
this feedback is embodied in the spiraling
repulsion of trajectories from the unstable convecting equilibria
at the center of each lobe or wing of the attractor before 
moving to the other lobe. 
Because the growth of oscillations is an important component to
the flow reversal process in both the CFD simulated thermosyphon
and EM system, we define here what is meant by oscillation
amplitude in each case.
In the CFD simulated thermosyphon,
it is the maximum distance of the system state from 
the nearest convecting equilibrium, where system state is understood
to mean the state of the entire temperature and velocity
flow fields in the CFD simulation.
When considering the DA-generated EM analysis, 
the $k$th $x_1$ oscillation amplitude
$x_1^\mathrm{max}$ is defined as the maximum amplitude
\begin{equation}
  \label{eq:x1oscamp}
  x_1^\mathrm{max} = \max_{t \in \mathcal{T}} | x_1(t) |
\end{equation}
where 
$\mathcal{T}=[t_0^{(k)}, t_f^{(k)}]$ is the time interval
of the $k$th oscillation.

\subsection{Flow reversal forecasting methods}
\label{sec:revoccur.tests}

Three separate tests were developed to predict, at each 
assimilation step, whether a 
flow reversal would occur within the next oscillation period
(approximately 11 min), here taken to be within the next
20 DA cycles. See Sec.~\ref{sec:revoccur.skill1} and  
Appendix S3 in the Supporting Information for a description of
how the tunable parameters were chosen.

\subsubsection{Lead forecast}

The simplest test forecasts a flow reversal
whenever the background forecast changes rotational state.
Note that to forecast a flow reversal occurring in the
future, the background forecast started from the most recent
analysis IC provides our only information
about the system's future state.
Ignoring the three-dimensional nature of the state space,
a flow reversal is forecast whenever $x_1$ crosses through
zero. Note also that the forecast is unable to predict flow reversals that occur
beyond the lead time, and that lead forecast quality quickly degrades as
the lead time is increased. We impose a
limit on the number of assimilation cycles to look ahead,
$\lambda_\text{lead} = 7$,
so that the algorithm does not
trust forecasts too far in advance.

\subsubsection{Bred vectors} 

An ensemble of perturbed states forming a small ball around the analysis
can be used to represent uncertainty in the IC. A nonlinear system will
dynamically stretch and shrink such a ball around its trajectory as it
moves through the attractor \citep{danforth2006}. Small perturbations
to points on a trajectory are integrated forward in time, and the
differences between perturbed and unperturbed solutions are called
bred vectors (BVs). Here, the rescaling amplitude is 0.001 and the
integration time coincides with the 30 s assimilation window.


The average BV growth rate is a useful measure of local instabilities
\citep{hoffman2009a}. 
\cite{evans2004}, studying
perfect-model forecasting of the Lorenz system, set a BV growth rate
threshold which accounted for 91.4\% of the observed flow
reversals (hit rate). Our BV test simply forecasts a flow reversal
whenever the average BV growth rate over the previous assimilation window
exceeds a threshold, $\rho_\text{BV} = 0.6786$.

\subsubsection{Correlation}
\begin{figure}
  \includegraphics[width=\linewidth]{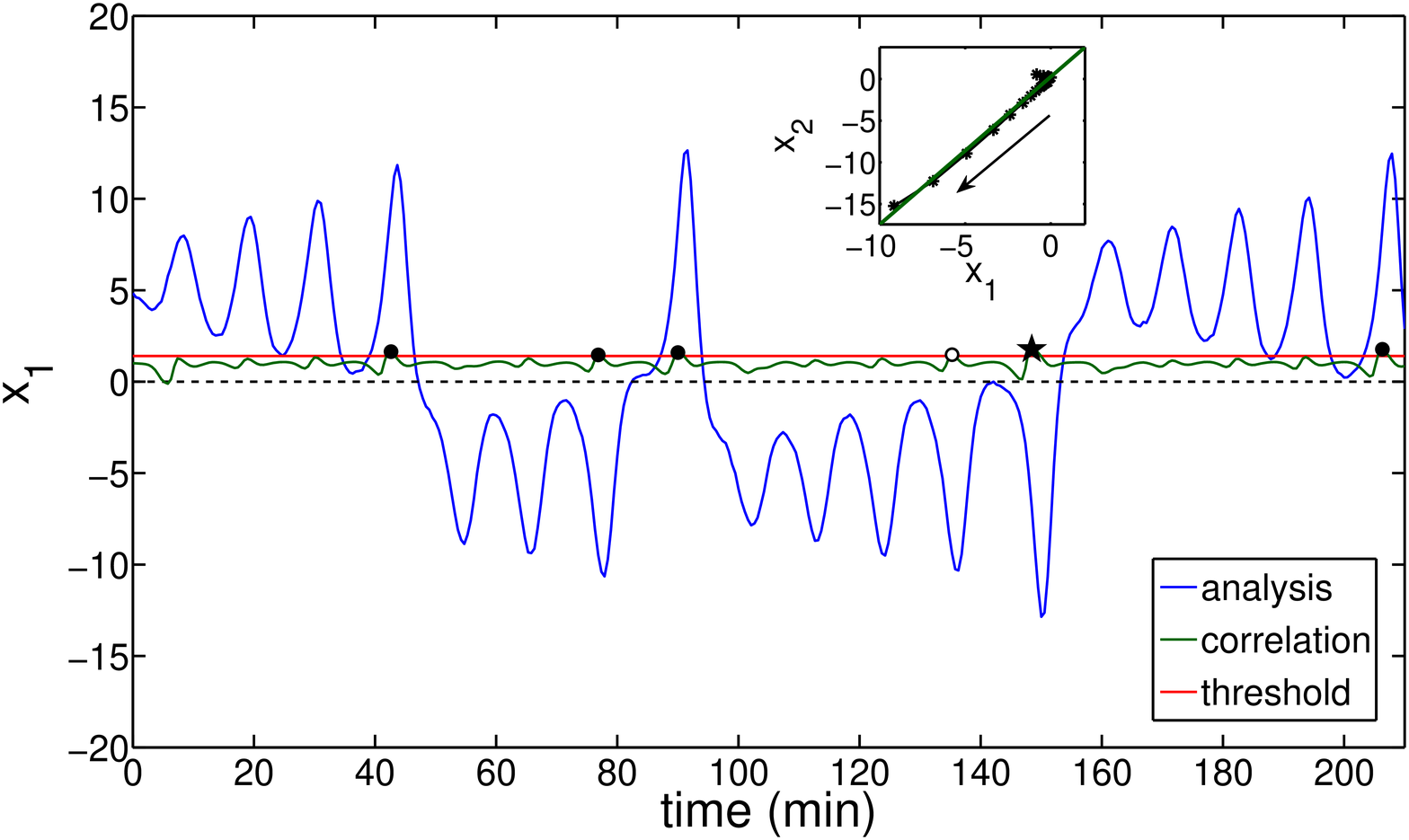}
  \caption{Correlation test: Whenever the slope (green) of the 
    $[x_2, x_1]^T$ correlation
    exceeds a threshold (red), a flow reversal is forecast. 
    Correct positive predictions are
    shown as filled circles and false positives as open circles. The
    starred point corresponds to the inset, which shows how correlation
    is computed as the slope of the least squares fit (green line) of
    previous analysis points, and the arrow shows 
    the direction of the trajectory.
    Note that each analysis cycle where the correlation fails
    to exceed the threshold counts as a correct negative
    forecast (not shown). There are no false negatives,
    i.e. misses, in this timeseries. 
    Here $\rho_\mathrm{corr}$ and $\lambda_\mathrm{corr}$
    are the same as for Tab.~\ref{tbl:con}.
  }
  \label{fig:corr}
\end{figure}

The final test uses the fact that flow reversals
are suspected to be caused by 
out of phase viscous and thermal dissipation. Since the
friction term grows with fluid velocity $\propto x_1$
and the thermal dissipation grows with
the size of the temperature anomaly, related to $x_2$, we examined the
correlation between those two variables over a tunable number of
previous analysis cycles. Specifically, when the slope of the
least-squares linear fit of $\lambda_\text{corr} = 18$ 
previous analysis points
$[ x_2(t-i), x_1(t-i) ]^T$
for 
$i=0,1,\ldots, ( \lambda_{\mathrm{corr}} - 1 )$
exceeds a threshold $\rho_{\mathrm{corr}} = 1.42$, 
a flow reversal is forecast. See Fig.~\ref{fig:corr} for an 
illustration of this process.
Interestingly, increasing autocorrelation of the state
seems to be a universal property of many systems in advance of critical
transitions \citep{scheffer2009, cotilla-sanchez2012a}.

\subsection{Forecasting residency times in the new rotational state}
\label{sec:revoccur.dur}

\begin{figure*}
  \centering
  \includegraphics[width=\linewidth]{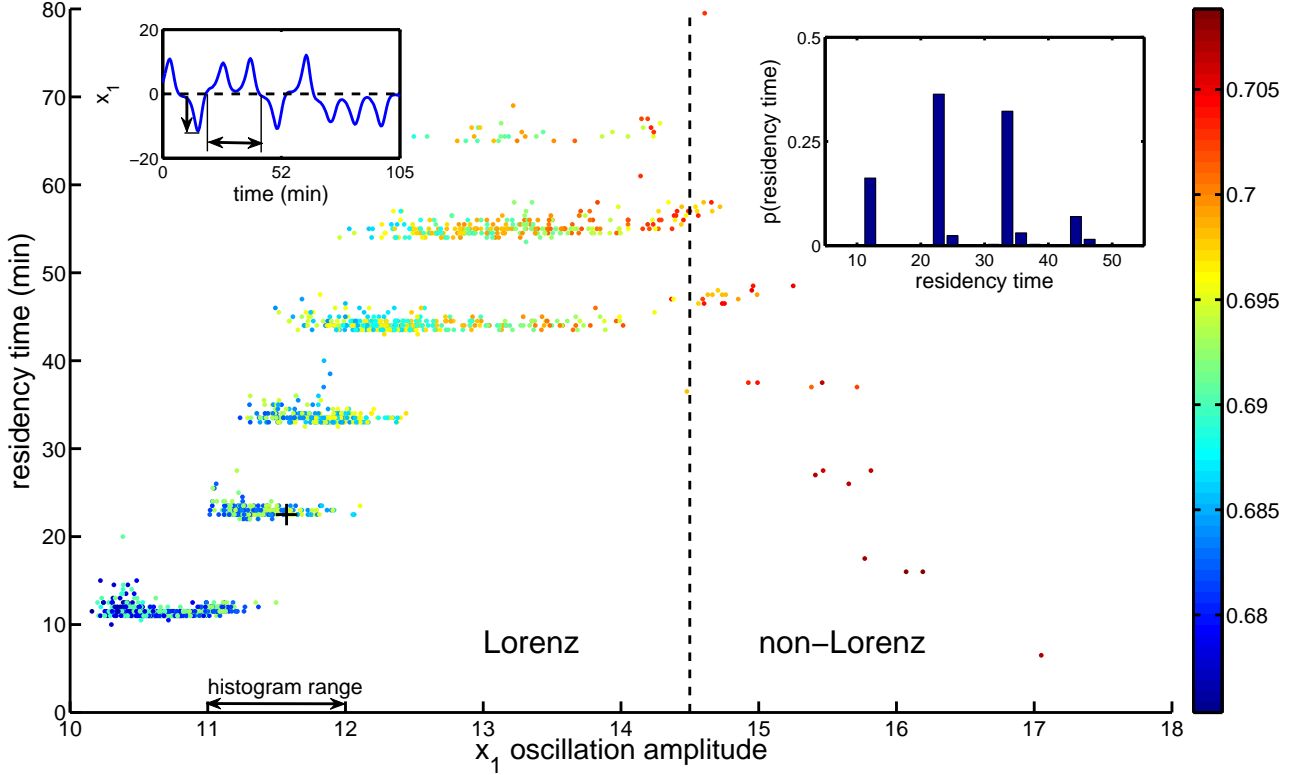}
  \caption{Residency time in the new rotational state is
    plotted versus the amplitude of $x_1$
    (proportional to the mass flow rate) at the last extremum before
    flow reversal. This amplitude is calculated from the EM
    model EKF analysis of thermosyphon observations, using a 30 s assimilation
    window. 
    This figure contains
    over 39 days of simulated flow and 1796 flow reversals. Points are
    colored by the average BV growth rate over the preceding assimilation
    window, showing a BV growth rate gradient that increases in the
    positive direction along both axes. 
    Inset, upper left: A timeseries corresponding to a single point
    in the scatter plot, marked with a black cross. The $x_1$ oscillation
    amplitude and subsequent residency time are denoted by the vertical
    and horizontal arrows, respectively.
    Inset, upper right: 
    A histogram
    showing the likelihood of residency times following
    an $x_1^\text{max} \in (11, 12)$. This is the interval
    that we would consider for an $x_1$ oscillation amplitude of 10.5
    preceding flow reversal.
    The most likely residency time is about 23
    min or 2 oscillations, the middle ``step'' for the histogram range.
    \label{fig:EMFregchange} 
  }
\end{figure*}

We found that the analysis' $x_1$ oscillation amplitude preceding
each flow reversal is correlated
with the duration of the following rotational state,
shown in Fig.~\ref{fig:EMFregchange}. 
We refer to these durations between flow reversals 
as residency times.
Residency times are observed at discernible
``steps'' corresponding to integer numbers of oscillations. 
This correlation makes the $x_1$ oscillation amplitude a
plausible predictor for residency time in the new rotational state.

Furthermore, the average BV growth
rate measured over the assimilation window preceding that extremum
follows a clear gradient in Fig.~\ref{fig:EMFregchange},
the growth rate increasing with
oscillation amplitude. The BV growth rate gradient implies that more
unstable system states precede longer residency times in the next
rotational state. 
Outliers with 
$x_1^\mathrm{max} \gtrapprox 14.5$
result in shorter residency times than expected from making similar plots
to Fig.~\ref{fig:EMFregchange} 
for the pure Lorenz and EM systems (not shown). In the Lorenz and
EM systems, the steps continue to move upwards with $x_1$ oscillation
amplitude.
The discrepancy is due to the non-Lorenz behavior that was
mentioned at the end of Sec.~\ref{sec:da.results}

Our residency time prediction algorithm proceeds as follows.
When a flow reversal is forecast by one of the methods
described in Sec.~\ref{sec:revoccur.tests}, the algorithm
first calculates $x_1^\text{max}$
as defined by Eqn.~\eqref{eq:x1oscamp}
for the presently occurring oscillation. 
The algorithm uses only the
analysis and lead forecast data available at the time 
the flow reversal test is triggered when estimating $x_1^\text{max}$.
The algorithm then
examines the residency times of all flow reversals which followed an
$x_1^\text{max}$ in the interval $(x_1^\text{max}-0.5, x_1^\text{max}+0.5)$.
These ordered pairs of amplitudes and residency times are drawn from the 
training timeseries.
From the relative abundance of residency times in this sample, 
we assign a probability to the number of flow oscillations in the
forthcoming rotational state. 
(See the inset histogram
in Fig.~\ref{fig:EMFregchange}.)
The categories are restricted to
1--6 oscillations (a duration of 7 oscillations,
shown in Fig.~\ref{fig:EMFregchange},
is observed exactly once in the training timeseries, so
it was considered too rare an event to merit a category).
The typical residency times corresponding to 1, 2, 3, 4, 5, and 6 oscillations
are taken to be 11.48, 23.09, 33.72, 44.38, 55.11, and 66.08 minutes, 
respectively; the oscillation category associated with a given
residency time is taken to be that with the closest time in 
this list.
This algorithm generates a probabilistic forecast from the relative 
abundance of points in each oscillation category.
An example output would be
20\%, 40\%, 30\%, and 10\% chance of 1, 2, 3, and 4 oscillations
in the next rotational state and zero probability of 5 or 6 oscillations.

\subsection{New details regarding the flow reversal mechanism}
\label{sec:revoccur.new}

\begin{figure*}



  \includegraphics[width=.75\linewidth]{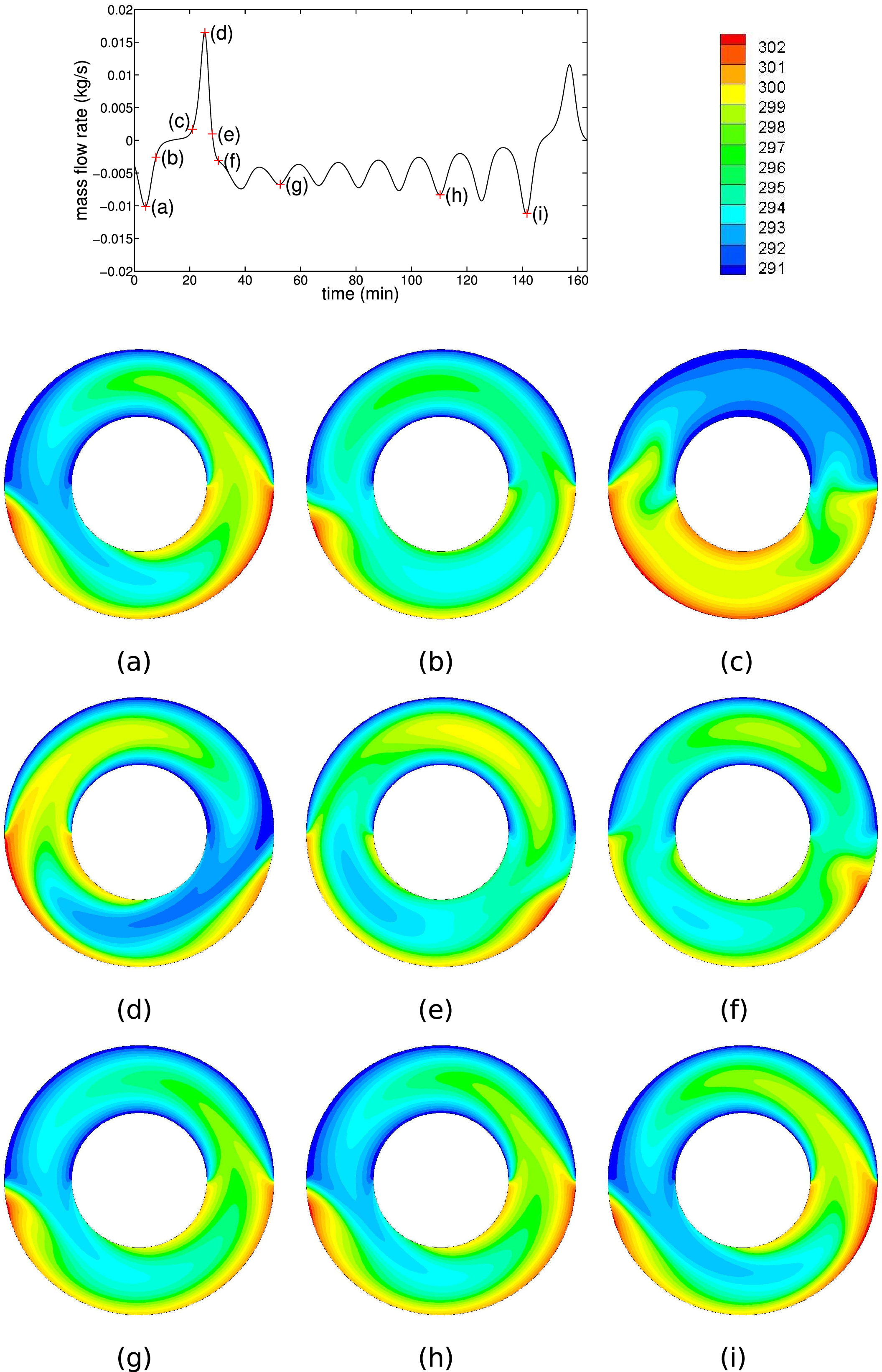}

  \caption{
    Temperature profiles before and after two flow reversals
    for the chaotic regime in units of Kelvin. 
    For ease of visualization, the radius ratio is reduced to 3. 
    $\text{Ra}=1.8\times 10^4$.
    The amplitude of
    the oscillation shown in (a) is
    relatively weak; the temperature profile quickly approaches (b) that of
    conduction (c) where it heats up significantly before the reversal. 
    The extreme instability of the conducting state,
    near time 20 min,
    produces a large oscillation in
    the CW direction (d), immediately causing another flow reversal
    back to CCW. As the system enters the new rotational state,
    remnant heat stabilizes the flow [contrast (e-f) with (b-c)],
    necessitating a longer residency 
    in the new rotational state
    while the instability grows (g-i)     
    (note that for the radius ratio of 24
    no more than 7 oscillations are ever observed). 
    \label{fig:regchangeprofiles}
  }
\end{figure*}

Not all flow reversals occur when the system reaches the same
flow oscillation amplitude, 
nor do all rotational states last the same amount of
time. During a flow reversal, the fluid motion stalls after hot fluid
extends across the entire heating section into the cooling section
(see Sec.~\ref{sec:revoccur.trad} and Fig.~\ref{fig:regchangeprofiles}). 
The magnitudes of this hot ``tongue'' and, likewise, the opposite cold
tongue affect the stability of the system as it reverses. 
If the oscillation is small, it will mostly dissipate before
the new rotational state is entered, bringing the temperature profile close to
that of conduction. This is a highly unstable equilibrium, since the
vertical temperature gradient builds until the fluid in the bottom is
much hotter than the fluid above (illustrated in the analysis in
Fig. \ref{fig:perfectstorm}). When the fluid begins to rotate,
it accelerates rapidly. The large amount of heat carried by the fluid
brings the system state far from the convecting equilibrium. If the
oscillation is large (corresponding to a large deviation from
convecting equilibrium in temperature and velocity), remnant warm and
cool areas will be present in the top and bottom sections of the loop,
respectively. These stabilize the new rotational state near its convective
equilibrium. The resulting duration is longer since the instability
requires more time to grow before causing the next reversal. These two
situations are illustrated in 
Fig.~\ref{fig:regchangeprofiles}
and explain the trend in the Lorenz region of 
Fig.~\ref{fig:EMFregchange}.
Animations of the simulated
temperature field during flow reversal are consistent with this
explanation 
\footnote{A movie similar to the case shown in
  Fig. \ref{fig:regchangeprofiles} is provided online at
  \url{http://www.uvm.edu/\~kharris/thermosyphon/T-Ra-18000-new.mp4}.}.

We believe that the behavior in the extremely large oscillation,
non-Lorenz region,
where $x_1^\mathrm{max} \gtrapprox 14.5$ and shown in 
Fig.~\ref{fig:EMFregchange}, is caused by excessive
remnant thermal energy after flow reversal. Although the
temperature distribution present after a flow reversal is
configured in a way that stabilizes the flow in the new rotational
state, the very large magnitude of the temperature field is
a competing, destabilizing factor that dominates as 
$x_1^\mathrm{max}$ increases into the non-Lorenz region.
This leads to shorter durations in the new rotational state before
a second flow reversal occurs.

\subsection{Flow reversal forecast skill}
\label{sec:revoccur.skill1}

The results of the three tests are presented in Tab.~\ref{tbl:con}
as two-by-two contingency tables. Shown in Tab.~\ref{tbl:skill}
are the threat score
(TS), false alarm ratio (FAR), and probability of
detection (POD) \citep{wilks1995}. 
Given a non-probabilistic yes/no forecast with $a$ hits,
$b$ false alarms, $c$ misses, and $d$ correct negatives for a total
of $n$ events, these are defined as 
TS=$a/(a+b+c)$, FAR=$b/(a+b)$, POD=$a/(a+c)$.
Because flow reversals are relatively rare events,
the hit rate $(a + d)/n$ would be dominated
by correct negatives. Instead, TS is chosen as an
appropriate overall performance metric since it disregards
these frequent negative events
and takes into account both false alarms and misses.

There are trade-offs among the various skill scores for each flow
reversal test. Tuning the reversal tests then amounts to
multiobjective optimization,
attempting to
maximize TS, RPS-avg, and RPS-med (the skill scores used
for residency time forecasts, defined in Sec.~\ref{sec:revoccur.skill2}),
minimize FAR, and maintain POD above 95\%.
The goal was to tune each method to all-around good performance, for both 
reversal occurrence and residency time forecasts.
To guide the process,
plots of the skill scores were made for different tuning
parameters, but the final tuning was performed ad hoc. 
In Appendix S3 in the Supporting Information, Fig.~S3-2
shows one of these tuning experiments with the final
parameters chosen appearing in the center of each subfigure.

Considering TS alone, the lead forecast performed best, followed by
the correlation test, with the BV test performing poorest. The
BV test also had a very high FAR, leading us to conclude, in contrast
to the results of \cite{evans2004} for a perfect model experiment, 
that BV growth rate is a poor overall predictor
of flow reversals in a realistic thermosyphon. 
On the other hand, the correlation test had the
lowest FAR while maintaining a high TS, but this comes at the price of
more misses, resulting in a lower POD. 
The reasonable performance of
the correlation test in all areas lends circumstantial 
evidence to the claim
that out of phase dissipations are indeed the cause of flow reversals.

The flow reversal occurrence tests are triggered in different
situations, leading to variation in how far in advance flow reversals
are detected, the ``warning time''.
Warning times were only computed for hits, 
i.e. forecast flow reversals which were observed to occur.
The lead, BV, and correlation tests had average warning times of 175,
217, and 304 s respectively. Histograms of these warning times are
presented in Appendix S3 in the Supporting Information.

\begin{table}
  {\footnotesize
    (a) Lead forecast, n=175592
  }\\[1ex]
  \begin{center}
    \begin{tabular}{cc|c|c|}
      \cline{3-4} & &
      \multicolumn{2}{|c|}{Observed} \\ 
      \cline{3-4} & & Yes & No \\ 
      \hline
      \multicolumn{1}{|c|}{\multirow{2}{*}{Fcast}} & Yes & 4472 & 744 \\
      \cline{2-4} \multicolumn{1}{|c|}{} & No & 13 & 170363 \\
      \hline
    \end{tabular}\\[1ex]
  \end{center}
  {\footnotesize
    (b) BV test, n=121120
  }\\[1ex]
  \begin{center}
    \begin{tabular}{cc|c|c|}
      \cline{3-4} & &
      \multicolumn{2}{|c|}{Observed} \\
      \cline{3-4} & & Yes & No \\
      \hline
      \multicolumn{1}{|c|}{\multirow{2}{*}{Fcast}} & Yes & 4383 & 3203 \\
      \cline{2-4} \multicolumn{1}{|c|}{} & No & 102 & 121258 \\
      \hline
    \end{tabular} \\[1ex]
  \end{center}
  {\footnotesize
    (c) Correlation test, n=174925
  }\\[1ex]
  \begin{center}
    \begin{tabular}{cc|c|c|} \cline{3-4} & &
      \multicolumn{2}{|c|}{Observed} \\ 
      \cline{3-4} & & Yes & No \\ 
      \hline
      \multicolumn{1}{|c|}{\multirow{2}{*}{Fcast}} & Yes & 3540 & 239 \\
      \cline{2-4} \multicolumn{1}{|c|}{} & No & 945 & 170201 \\ 
      \hline
    \end{tabular} 
  \end{center}
  \vspace{1em}
  \caption{Contingency tables for the three flow reversal tests. 
    The parameters used were:
    $\rho_\mathrm{corr} = 1.42, \lambda_\mathrm{corr} = 18$,
    $\rho_\mathrm{BV} = 0.6786$, and $\lambda_\mathrm{lead} = 7$.
    For a detailed description of the tuning, see the text and 
    Fig.~S3-2 in Appendix S3 in the Supporting Information.
    \label{tbl:con}}
\end{table}

\begin{table}
  \centering
  \begin{tabular}{r||c|c|c|c|c}
    Method & TS & FAR & POD & RPS-avg & RPS-med \\
    \hline
    lead & 86 & 14 & $>$99 & 71 & 87 \\
    bred vector   & 57 & 42 & 98 & 67 & 86 \\
    correlation & 75 & 6  & 79    & 58 & 74 
  \end{tabular}\\[1ex]
  \caption{Skill scores for flow reversal categorical forecasts (TS, FAR, POD)
    and residency time probabilistic forecasts (RPS-avg, RPS-med) 
    for the three tests. Note that
    for TS and POD a perfect score is 100\%, while a perfect FAR is 0\%;
    RPS scores should be interpreted as a percent
    improvement over climatology, so that any RPS $>$ 0 is an improvement.}
  \label{tbl:skill}
\end{table}

\subsection{Residency time forecast skill}
\label{sec:revoccur.skill2}

A n\"aive way of forecasting residency times would assign each possible
outcome a probability equal to that measured from the climatology. 
In our case, this would amount to using
the marginal distribution of oscillation occurance.
However, our method also takes into account the $x_1^\text{max}$ before
the flow reversal (i.e. the joint distribution of events by 
oscillation occurance and $x_1^\text{max}$), 
which we have shown contains important information about
the number of oscillations that the system will undergo in the new 
rotational state.
So, we compare our method to climatology using
a ranked probability skill score (RPS, see \citealp{wilks1995}).
This is only computed in the case of hits.
We actually computed two variants, by taking either the mean (RPS-avg)
or median (RPS-med) of the ranked probability scores for each reversal
event when computing the skill.
The results are illustrated in Tab.~\ref{tbl:skill}. The lead forecast test
performs best, followed by the BV test and the correlation test.
Unsurprisingly, the flow reversal tests with smaller warning times
performed better when making residency time forecasts. Because there is
more information about the system state immediately preceding a
flow reversal if the warning time is small, the residency time
forecast is better informed. 
The magnitude of the improvement over climatology is large for
all methods. The RPS scores in Tab.~\ref{tbl:skill} are similar to
or better than those for probabilistic forecasts in NWP 
\citep{tippett2008a, doblas-reyes2000a}.

\section{Conclusion}
\label{sec:conclusion}

DA was shown to be an effective way of coupling a simplified
model to CFD simulations of the thermosyphon. Although
background forecast errors were always larger than observational noise, 
climatically scaled background error was small for reasonable
assimilation windows. Proper tuning of multiplicative and
additive inflation factors was essential for avoiding filter
divergence and achieving low forecast error. 
All of the DA methods used in this study accurately
capture the behavior of the thermosyphon with short assimilation windows.
Among the DA methods, the ensemble methods show advantages 
over 3D-Var and EKF with longer assimilation windows, when nonlinear
error growth becomes important.
With frequent analysis update, DA can reveal non-Lorenz behavior in the
thermosyphon even with the EM (Lorenz-like) model.

Three different predictors
of flow reversals were proposed and tested with reasonable success.
In comparison with the two rules in \cite{evans2004} 
for predicting the behavior of the Lorenz trajectory, 
the BV growth rate is a useful measure of the EM model's 
dynamical instabilities, but it does not perform well on its
own as a predictor of flow reversals.
Finally, the amplitude of the final oscillation in
the current rotational state was found to be correlated
with the residency time in the following rotational state, 
and we provide a physical explanation
for this phenomenon, elaborating on the details of flow reversals.
Oscillation amplitudes were then used to create
probabilistic forecasts of those residency times with significant
improvements over climatology.

A laboratory thermosyphon device is in construction. The next stage
of this research will apply similar methods to forecasting the system state,
flow reversals, and residency times using
3D numerical flow simulations. 
Spatially-aware DA techniques, such as the 
Local Ensemble Transform Kalman Filter \citep{kalnay2007a,hunt2007},
could be applied to finite-volume or finite-element models to
study the spatial structure of the fluid flow and error growth.
These imperfect model experiments could be used
to compare the relative performance of other DA algorithms 
(4D-Var, \citealp{kalnay2007a}),
synchronization approaches (adaptive nudging, see
\citealp{yang2006a}), 
and empirical correction techniques 
\citep{danforth2007a, li2009a,allgaier2012a}.

Although the thermosyphon is far from representing anything as complex
and vast as Earth's weather and climate, there are
characteristics our toy climate shares with global atmospheric
models.
Sophisticated atmospheric models are, at best, 
only an approximate representation of the
numerous processes that govern the Earth's climate. 
Global weather models and
the EM model both parameterize fine-scale processes that interact
nonlinearly to determine
large-scale behavior. 
Clouds and precipitation are sub-grid scale 
processes in a global weather model, and the correlations
for the heat transfer and friction coefficients are parameterizations
of fluid behavior on a finer scale than can be dealt with
in the reduced model.
Cloud formation is only partly understood, and moist
convection is an area of active research where some models bear
similarities to the EM model considered here \citep{weidauer2011a}.

The methods we use to forecast the toy model are also similar to the
methods used for global geophysical systems. Both require state
estimation to find the IC from which to generate forecasts. Also, when
forecasts are made in either system, climatology and
dynamically accessible regimes are often more
important than specific behavior:
the occurrence of flow reversals for the thermosyphon; 
periodic behavior such as the El Ni\~{n}o Southern Oscillation,
and statistics such as globally and regionally-averaged temperatures and their
effects on rainfall, ice cover, etc. for climate. 
Each of these is
a statistic that must be post-processed from the model output. To
meet these global challenges, many tools are needed in the modeling
toolbox. 
These techniques
may also be useful for forecasting sociotechnological systems
which are
rapidly gaining importance as drivers of human behavior.
In this way, toy models can provide us with insights 
that are applicable to the important scientific problems of today.

\section*{Acknowledgments} 
We would like thank Dennis Clougherty, Peter
Dodds, Nicholas Allgaier, and Ross Lieb-Lappen
for comments and discussion and the three anonymous reviewers
for providing many comments and suggestions that strengthened the paper.
We also thank Shu-Chih Yang
for providing 3D-Var MATLAB code
that was used to prototype our own
experiments.
We also wish to acknowledge financial support from the Vermont Space Grant
Consortium, NASA EPSCoR, NSF-DMS Grant \#0940271, the
Mathematics \& Climate Research Network, 
and the Vermont Advanced Computing Center.


\appendix
\renewcommand{\theequation}{S\arabic{section}-\arabic{equation}} 
\renewcommand{\thesection}{S\arabic{section}}
\renewcommand{\thesubsection}{S\arabic{section}.\arabic{subsection}}
\renewcommand{\thefigure}{S\arabic{section}-\arabic{figure}}
\renewcommand{\thetable}{S\arabic{section}-\arabic{table}}
\setcounter{figure}{1}
\setcounter{table}{1}

\section{Ehrhard-M\"uller model}
\subsection{Derivation}

\label{sec:model.derivation}
Following the derivations of 
\cite{gorman1986} and \cite{ehrhard1990}, we
consider the forces acting upon a control volume of 
incompressible fluid in the loop. 
All fluid properties are cross-sectionally averaged, and the
radial components of velocity and heat conduction within the fluid are
neglected. The fluid velocity $u=u(t)$ is assumed to be constant at all
points. 
Applying Newton's
second law, the sum of all forces on the control volume must equal its
change in momentum:
\begin{subequations}
  \label{eq:EMforces}
  \begin{equation}
    \label{eq:Feqma} F_p + F_f + F_g = \rho \pi r^2 R d\phi
    \frac{du}{dt}
  \end{equation} 
  where
  \begin{align} F_p &= -\pi r^2 R \, d\phi \nabla p = -\pi r^2 d\phi
    \frac{\partial p}{\partial \phi} \label{eq:Fp}\\ 
    F_f &= -\rho \pi r^2
    R \, d\phi \, f_w \label{eq:Ff}\\ 
    F_g &= -\rho \pi r^2 R \, d\phi \, g
    \sin{\phi}\,.
    \label{eq:Fg}
  \end{align}
\end{subequations} 
The angular coordinate $\phi$ and loop dimensions $r$ and $R$ are defined
in Fig.~1 in the main text, $g$ is the acceleration of gravity,
$\rho$ is the fluid density, $u$ is velocity, and $p$ is pressure.
The total force in Eqn.~\eqref{eq:Feqma} is
comprised of the net pressure ($F_p$), friction from shear within the
fluid ($F_f$), and the force of gravity ($F_g$). The pressure term,
Eqn.~\eqref{eq:Fp}, is the volume times the pressure gradient. 
The friction term, Eqn.~\eqref{eq:Ff}, is written in this form in
order to simplify the analysis; all frictional effects are contained
in $f_w$ which will depend on fluid velocity, to be discussed later.

Before we write the momentum equation, it is convenient to apply the
Boussinesq approximation, which assumes that variations in fluid
density are linear with temperature. In other words, $\rho=\rho(T)
\approx \rho_0 (1-\gamma (T-T_0))$ where $\rho_0$ is the reference
density, $\gamma$ is the coefficient of volumetric thermal expansion,
and $T_0=\frac{1}{2}(T_h+T_c)$ is the reference temperature. The
Boussinesq approximation also states that the density variation is
insignificant except in terms multiplied by $g$. Thus, the density
$\rho$ is replaced by $\rho_0$ in all terms of
Eqn.~\eqref{eq:EMforces} except gravity, Eqn.~\eqref{eq:Fg}. Using the
Boussinesq approximation, gathering terms, and dividing out common
factors gives the momentum equation
\begin{equation}
  \rho_0 \frac{du}{dt} d\phi = 
  -d\phi
  \left(\frac{1}{R} \frac{\partial p}{\partial \phi} + \rho_0
    \left(1-\gamma(T-T_0)\right) g \sin \phi + f_w \right) \:.
\end{equation}
Integrating about the loop, the momentum equation is
simplified because $u$ and $f_w$ are independent of $\phi$ and other
terms drop out due to periodicity.
\begin{align} \rho_0 \frac{du}{dt} &= \frac{\rho_0 \gamma g}{2\pi}
  \int_0^{2\pi}\!\!  d\phi \; T \sin \phi -f_w
  \label{eq:momentum}
\end{align}

We now must account for the transfer of energy within the fluid, and
between the fluid and the wall. All modes of heat transfer are neglected
except convection, which is a valid approximation when $r \ll R$ 
\citep{welander1967,ehrhard1990}.
The energy rate of change ($D/Dt$ is the material derivative
with respect to time) in the control volume is 
\begin{align} \label{eq:deltae} \rho_0 \pi r^2 R \, d\phi \, c_p
  \frac{DT}{Dt} & \equiv \rho_0 \pi r^2 R \, d\phi \, c_p \left(
    \frac{\partial T}{\partial t} + \frac{u}{R}\frac{\partial T}{\partial
      \phi} \right)
\end{align} 
which must be equal to the heat transfer through the wall
\begin{align}\label{eq:deltaq} \Delta Q &= -\pi r^2 R \, d\phi \, h_w
  (T-T_w) \,,
\end{align} 
where $c_p$ is the specific heat of the fluid,
$h_w$ is the heat transfer
coefficient, which depends on velocity, and $T_w$ is the temperature
at the wall. Combining Eqns.~\eqref{eq:deltae} and~\eqref{eq:deltaq}
gives the energy equation
\begin{align} \label{eq:energy} \left( \frac{\partial T}{\partial t} +
    \frac{u}{R}\frac{\partial T}{\partial \phi} \right) &=
  -\frac{h_w}{\rho_0 c_p} (T-T_w) \,.
\end{align} Together, Eqns. \eqref{eq:momentum} and \eqref{eq:energy}
represent a simple model of the flow in the loop.

The transport coefficients $f_w$ and $h_w$ characterize the
interaction between the fluid and the wall. They are defined by the
constitutive relations \citep{ehrhard1990}
\begin{align} h_{w} &= h_{w_0} \left[ 1 + K h(|x_1|)
  \right] \label{eq:h_w}\\ f_{w} &= \frac{1}{2} \rho_0 f_{w_0} u \,
  , \label{eq:f_w}
\end{align} 
where $x_1 \propto u$ is the dimensionless velocity. The function
\begin{equation}
  \label{eq:H}
  h(x) = \left\{
    \begin{array}{lr}
      p(x)    & \text{when } x < 1 \\
      x^{1/3} & \text{when } x \geq 1
    \end{array}
  \right. 
\end{equation}
in Eqn.~\eqref{eq:h_w}
determines the velocity dependence of the heat transfer coefficient,
which varies as $u^{1/3}$ for moderate $u$ \citep{ehrhard1990}. We
introduce the fitting polynomial $p(x)=44/9 \,x^2 - 55/9 \,x^3 +
20/9\, x^4$ to ensure that $h_w$ is analytic at $x_1=0$. 
This piecewise definition causes $h(x)$ to vary as $p(x)$
for $x \leq 1$ and $x^{1/3}$ for $x>1$. Eqn.~\eqref{eq:f_w} gives the
frictional deceleration of the fluid when $|u|>0$, and the $\rho_0/2$ term 
is retained to simplify the final solution. Dimensionally, $f_w$ is an
acceleration (m/s$^{2}$) and $h_w$ is power per unit volume per unit
temperature (W/m$^{3}$K). These coefficients $h_{w_0}$, $f_{w_0}$, and $K$
must be estimated from experiments (e.g., 
\citealp{ehrhard1990,welander1967,gorman1986}) or from other
empirical means. In Sec.~\ref{sec:model.paramest},
we describe the empirical methods used for parameter estimation.

\cite{ehrhard1990} solved the system of two
coupled, partial differential equations (Eqns. \eqref{eq:momentum} and
\eqref{eq:energy}) by introducing an infinite Fourier series for
$T$. The essential dynamics can be captured by the lowest modes, i.e.,
\begin{align}\label{eq:T} 
  T(\phi,t) &= C_0(t)+S(t) \sin \phi + C(t) \cos \phi \,.
\end{align} 
Because this form of $T$ separates the variables $\phi$ and
$t$, the problem is transformed into a set of ordinary differential
equations. Substituting Eqn.~\eqref{eq:T} into
Eqn.~\eqref{eq:momentum} and integrating gives the equation of motion
for $u$. Similarly, Eqn.~\eqref{eq:energy} is integrated by $\oint
d\phi \sin \phi$ and $\oint d\phi \cos \phi$ to separate the two
temperature modes $S $ and $C$. The system is written in dimensionless
form
\begin{subequations} \label{eq:EM}
  \begin{align} \frac{dx_1}{dt'} &= \alpha \left(x_2 - x_1\right)\\
    \frac{dx_2}{dt'} &= \beta x_1 - x_2 \left( 1+K h(|x_1|) \right) - x_1
    x_3\\ \frac{dx_3}{dt'} &= x_1 x_2 - x_3 \left( 1+K h(|x_1|) \right)
  \end{align}
\end{subequations} 
where the following linear transformations have
been made to create dimensionless variables
\begin{equation} \label{eq:EMdim} \left.
    \begin{aligned} t' & = \frac{h_{w_0}}{\rho_0 c_p}t\\ x_1 & =
      \frac{\rho_0 c_p}{R h_{w_0}}u\\ x_2 & =\frac{1}{2} \frac{\rho_0 c_p
        \gamma g}{R h_{w_0} f_{w_0}} \Delta T_{3-9}\\ x_3 & =\frac{1}{2}
      \frac{\rho_0 c_p \gamma g}{R h_{w_0} f_{w_0}}
      \left(\frac{4}{\pi}\Delta T_w - \Delta T_{6-12} \right) \;.
    \end{aligned} \quad \right\}
\end{equation}

Physically, $x_1$ is proportional to the mean fluid velocity, $x_2$ to
the temperature difference across the convection cell or $\Delta
T_{3-9}$ (between 3 o'clock and 9 o'clock), and $x_3$ is proportional
to the deviation of the vertical temperature profile (characterized by
the temperature difference between 6 o'clock and 12 o'clock, $\Delta
T_{6-12}$) from the value it takes during conduction.

The parameter $\alpha = \frac{1}{2} \rho_0 c_p f_{w_0} / h_{w_0}$ is
comparable to the Prandtl number, the ratio of momentum diffusivity
and thermal diffusivity. Similar to the Rayleigh number, the heating
parameter
\begin{align} \beta &=\frac{2}{\pi} \frac{\rho_0 c_p \gamma g}{R
    h_{w_0} f_{w_0}} \Delta T_w
\end{align} determines the onset of convection as well as the
transition to chaotic behavior.

Although the previous derivation assumes a 3D geometry,
the CFD simulations described in Sec.~2.1 of the main text,
were performed in 2D.
A 2D geometry corresponds to infinite concentric cylinders as opposed to
the quasi-1D torus. Due to cross-sectional averaging, the EM equations of
motion \eqref{eq:EM} are the same in 2D or 3D; the change may be 
realized by letting $\pi r^2 \to 2r$ in 
Eqns.~\eqref{eq:EMforces}, \eqref{eq:deltae}, and \eqref{eq:deltaq}
and carrying out the rest of the derivation.
The only differences arise in the non-dimensional transformations
and parameters, which were empirically determined by a
multiple shooting algorithm (see Sec.~\ref{sec:model.paramest}).

\subsection{Parameter estimation} 

\label{sec:model.paramest} 

Before any forecasting, the parameters matching the EM model to
the thermosyphon simulation needed to be determined. 
\cite{ehrhard1990} used experimental measurements
to determine the correlation
coefficients for friction, $f_{w_0}$, and heat transfer, $h_{w_0}$ and
$K$. They achieved this by opening the loop at $\phi=\pi/2$ and
providing a developed flow with adjustable velocity.  By measuring the
pressure loss ($\propto f_w$) and 
heat transfer across the loop for a range of
velocities, they were able to find the correlation coefficients using
regression. We were unable to accomplish this with a CFD
simulation of an open-loop geometry.

Instead, parameter estimation was formulated as a multiple
shooting problem. Shooting methods minimize the error in 
an ODE trajectory relative to data by optimizing over
all possible initial conditions and parameter space.
Multiple shooting is a shooting method suitable for chaotic
ODEs \citep{baake1992a}.
It overcomes the sensitive dependence on initial conditions
by partitioning the data set and
solving the shooting problem on those subsets of the data, 
augmented by continuity conditions.
We used the nonlinear least square optimizer
\texttt{lsqnonlin} in \cite{matlab} to perform the minimization
and relaxed the continuity constraints.
The model parameters which were tuned were $\alpha$, $\beta$, and $K$. 
However, we also needed a way to determine the time and velocity
scales to convert the dimensionless variables $t'$ and $x_1$ to
their observed, dimensional values $t$ and $q$. 
These scales change as the other
parameters are varied, so these were incorporated
into the variables of the optimization.

\begin{figure}
    \centering
    \includegraphics[width=\linewidth]{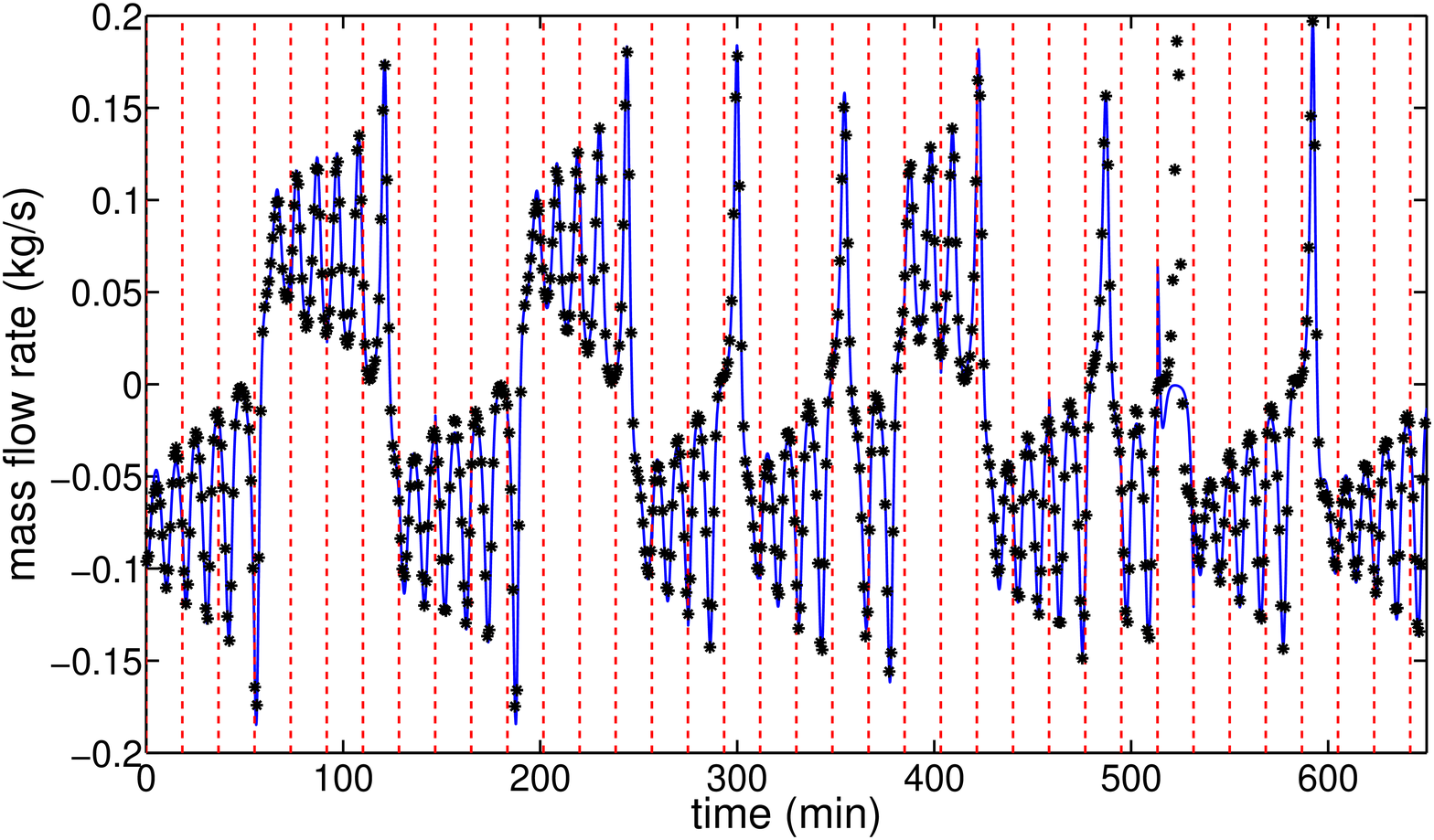}
    \caption{
      Results of the multiple shooting algorithm. 
      The entire training timeseries is shown.
      Starred points
      are the CFD thermosyphon mass flow rate data used for validation. 
      Also shown are the trajectories of
      model integrations over each shooting window, which are indicated
      by the dashed vertical lines. In one shooting window,
      near 515 min, the model does not match the data because the
      optimization has not found a good IC 
      (it is near the highly unstable conducting equilibrium)
      even though the parameters are acceptable.
    }
    \label{fig:multshoot}
\end{figure}

\begin{table}
  \centering
  \begin{tabular}{l|l}
    parameter & value \\
    \hline
    $\alpha$       & 7.99 \\
    $\beta$        & 27.3 \\
    $K$            & 0.148 \\
    $t$ scale (s)  & 631.6 \\
    $q$ scale (kg/s) & 0.0136
  \end{tabular}
  \caption{
    Final parameters used in the data assimilation scheme. The
    first three are the dimensionless parameters of the model, and
    the final two are used to rescale the dimensionless time and mass
    flow rate.}
  \label{tbl:params}
\end{table}

The results of the multiple shooting algorithm are shown in 
Fig.~\ref{fig:multshoot} and Table~\ref{tbl:params}.

\section{Data assimilation algorithms} \label{sec:ax:da}

\subsection{Kalman Filter (KF)}

The KF is well-known and widely used in linear DA and control
problems. Although the thermosyphon is highly nonlinear, the linear
update equations are similar to those of the nonlinear algorithms used
for this experiment. The KF attempts to assimilate observations and
forecasts for a process of the form
\begin{align} \xvec_{k}^t &= \mathbf{W} \xvec_{k-1}^t \,.
\end{align} In this case, $\xvec^t$ is the true state, which advances
in time according to the linear process $\mathbf{W}$, which is unknown
but approximated by the model $\mathbf{M}$. Subscripts index the time
step. Using the model, the analysis from the previous time step is
integrated to generate the background forecast for the current time
step
\begin{equation} \label{eq:analysisfcast} \mathbf{x}_{k}^b =
  \mathbf{M} \mathbf{x}_{k-1}^a
\end{equation} where $\mathbf{M}$ is the linear model, $\xvec^a$ is
the old analysis, and $\xvec^b$ is the background. Because
$\mathbf{M}$ is only an approximation of $\mathbf{W}$, a perfect
initial condition will not lead to a perfect forecast, so
\begin{equation} \xvec_k^t = \mathbf{M}\xvec_{k-1}^t + \epsilon_k^q
\end{equation} where the model errors $\epsilon^q$ have covariance
$\mathbf{Q}$ (usually assumed to be constant in time) and are written
on the right hand side for convenience. When deemed unnecessary, time
subscripts are left out.

Given an observation $\yvec$ and background forecast $\xvec^b$, the KF
finds the optimal way to combine them into the \textit{analysis}
$\xvec^a$, the best guess of the current state. This becomes the IC
when forecasting with the model, Eqn.~\eqref{eq:analysisfcast}. In an
operational context, we usually cannot observe every state
variable. If $\xvec \in \mathbb{R}^N$ and $\yvec \in \mathbb{R}^M$,
then $M < N$ (in NWP $M \ll N$), so we define the observation operator
$\mathcal{H}: \mathbb{R}^N \to \mathbb{R}^M$ that takes the background
forecast from the model state space into the observation space. This
serves two purposes: first, it avoids extrapolation of observations to
gridpoints in state space; and second, it enables us to interpret our
forecasts by comparing them directly to observations. For the
thermosyphon, $\mathcal{H}$ is linear, so we write it as $\mathbf{H}$,
but this is usually not the case for the observations in NWP,
e.g., satellite radiances and radar reflectivities.

The complete application of the KF consists of a forecast step
\begin{subequations}
  \label{eq:KFfcast}
  \begin{align} \xvec^b_{k} &= \mathbf{M}
    \xvec^a_{k-1} \label{eq:KFfcaststate}\\ \mathbf{B}_{k} &= \mathbf{M}
    \mathbf{A}_{k-1} \mathbf{M}^T + \mathbf{Q} \label{eq:KFfcastcov}
  \end{align}
\end{subequations} and an analysis step
\begin{subequations}
  \label{eq:KFanalysis}
  \begin{align} \xvec^a_{k} &= \mathcal \xvec^b_{k} + \mathbf{K}_{k}
    (\yvec_{k} - \mathbf{H} \xvec_{k}^b) \label{eq:KFanalysisvec}\\
    \mathbf{A}_{k} &=
    (1-\mathbf{K}_{k}\mathbf{H})\mathbf{B}_{k} \label{eq:KFanaylsiscov}
  \end{align}
\end{subequations} with the Kalman gain $\mathbf{K}_k$ given by
\begin{align} \label{eq:KFgain} \mathbf{K}_k=\mathbf{B_\mathit{k}
    H}^T(\mathbf{H B_\mathit{k} H}^T+\mathbf{R})^{-1} \,.
\end{align}

The forecast equations create the background forecast and update the
background error covariance. The new background error covariance is
the old analysis error integrated forward plus the model error
$\mathbf{Q}$. In the analysis step, this background forecast is
incremented by the gain times the innovation ($\yvec - \mathbf{H}
\xvec^b$) to produce the analysis. The difference between the analysis
and the background is referred to as the \textit{analysis increment};
statistical properties of these increments can be used to reduce model
error \citep{danforth2007a, danforth2008a, danforth2008b}. The new
analysis error is equal to the background error reduced by a factor of
$(1-\mathbf{K}\mathbf{H})$. By finding the analysis, the filter has
revealed the best possible starting point for the next background
forecast. In fact, if the system is linear, the KF is the optimal
algorithm for state-estimation.

\subsection{Variational Filtering (3D-Var)}

Rather than minimize the analysis error variance, the analysis
equations can also be derived by finding the analysis state $\xvec^a$
that minimizes the quadratic scalar cost function $2C(\xvec) =
(\xvec-\xvec^b)^T \mathbf{B}^{-1}
(\xvec-\xvec^b)+(\yvec-\mathcal{H}\xvec)^T
\mathbf{R}^{-1}(\yvec-\mathcal{H}\xvec)$.  The cost $C(\xvec)$ has its
minimum at $\xvec=\xvec^a$, where $\xvec^a$ is given by
Eqn.~\eqref{eq:KFanalysis}. This is called the 3D variational (3D-Var)
method since the minimization for NWP is with respect to a state
vector embedded in a three-dimensional field (latitude, longitude, and
height).

Formally, both 3D-Var and the KF yield the same solution
\citep{kalnay2002}. However, in this case the control variable is the
analysis, while in the KF the control variable is the weight matrix
itself. In operational NWP, where the dimension of the state space $N$
is of $\mathcal{O}(10^9)$, the numerical implementations of 3D-Var and
the nonlinear KF are drastically different. Because 3D-Var assumes the
background error $\mathbf{B}$ is fixed in time, the Kalman gain
$\mathbf{K}$ needs to be calculated only once. The calculation of
$\mathbf{K}$ is the most computationally prohibitive part of DA
because it requires solving a linear system in $N$ variables. A
constant $\mathbf{K}$ thus makes the algorithm computationally simple;
the most difficult part of implementing 3D-Var is finding the optimal
$\mathbf{B}$.

However, a static $\mathbf{B}$ is not realistic. From a dynamical
systems standpoint, uncertainty is closely related to stability, which
is clearly dependent on the system state. In the thermosyphon, the
true background error is typically smaller when the system state is
near the unstable convecting equilibria than when the state is near
the more unstable conducting equilibrium. Because 3D-Var is
computationally cheap, the National Centers for Environmental
Prediction (NCEP) employ it to estimate ICs for the National Weather
Service 14-day global forecasts. However, it cannot detect so-called
``errors of the day'', state-dependent forecast errors which grow
quickly but are not represented in the 3D-Var background error
covariance matrix \citep{kalnay2002,li2009a}.

In our implementation, the 3D-Var background error covariance
$\mathbf{B}_\text{3D-Var}$ was calculated iteratively,
using a techniques similar to that described in \cite{yang2006}.
We did this by calculating a time average of
the outer product of analysis increments
\begin{equation}
  \label{eq:B3dvar}
  \mathbf{B}_\text{3D-Var} = 
  \langle (\xvec^b - \xvec^a)(\xvec^b - \xvec^a)^T \rangle ,
\end{equation}
disregarding the initial 500 assimilation cycles and iterating the
process until convergence. 
During this, forecast errors were observed to decrease and stabilize.
This $\mathbf{B}_\text{3D-Var}$ was first computed for 
the 30 s assimilation window. It was stored and then used to bootstrap
the iterative procedure for the 60 s assimilation window, which
was stored and fed into the calculation for the 90 s 
assimilation window, etc.

\subsection{Extended Kalman Filter (EKF)}

The EKF is essentially the KF applied to a nonlinear model. Given a
nonlinear model $\mathcal{M}$, the error covariances are updated by
the {\it linear tangent model} $ \mathbf{M} \equiv \left.\partial
  \mathcal{M}/\partial \xvec\right|_{\xvec=\xvec^b} $ which takes the
place of $\mathbf{M}$ in Eqn.~\eqref{eq:KFfcastcov}. This model
propagates small perturbations around the trajectory $\xvec^b$ forward
in time. To operate on the matrix $\mathbf{A}$ with the linear tangent
model, first take the Jacobian of $F$ (the right hand side of the
nonlinear differential equation $\dot{\xvec}=F(\xvec)$ which describes
the model $\mathcal{M}$) and evaluate it at the background point
$\xvec^b$; call this matrix $\mathbf{J}$. Each column $\mathbf{a}_i$
of $\mathbf{A}$, which can be thought of as an error perturbation to
the analysis state, is then integrated forward in time according to
the linear ODE $\dot{\mathbf{a}_i}=\mathbf{J} \, \mathbf{a}_i$.

Also note that if the observation operator $\mathcal{H}$ is nonlinear,
it is replaced by a similar linear tangent model $\mathbf{H}$ in the
matrix equations \eqref{eq:KFanalysis} and \eqref{eq:KFgain}. The
transpose of these matrix functions are called {\it adjoint models},
which are used in sensitivity analysis of the state to perturbations.

To propagate the background covariance without the explicit adjoint
model, as Eqn.~\eqref{eq:KFfcastcov} would require, $\mathbf{B}$ was
first decomposed with the Cholesky factorization \citep{golub1996} into
the product of a lower and upper diagonal matrix before its columns
were integrated forward with the linear tangent model $\mathbf{M}$.
\begin{align} \mathbf{B}_{k-1} &= \mathbf{L}_{k-1}
  \mathbf{L}_{k-1}^{T}\\ \mathbf{T}_k &= \mathbf{M}_{k-1}
  \mathbf{L}_{k-1} \\ \mathbf{A}_k &= \mathbf{T}_k \mathbf{T}_k^{T} +
  \mathbf{Q}
\end{align} This guarantees symmetry for the new analysis error
covariance $\mathbf{A}$.

Some modifications to the EKF algorithm are necessary to prevent
filter divergence.
A multiplicative inflation factor 
\begin{equation}
  \label{eq:inflmultekf}
  \mathbf{B} \leftarrow (1+\Delta) \mathbf{B}
\end{equation}
was applied to the
background covariance matrix after the model integration
and before the analysis step.  
We also performed additive inflation,
following \cite{yang2006}. Random numbers uniformly
distributed between 0 and $\mu$ were added to the diagonal elements of
$\mathbf{A}$ after performing the analysis and before the
next forecast step, i.e.
\begin{equation}
  \label{eq:infladdekf}
  \mathbf{A} \leftarrow \mathbf{A} + \mu \; \mathrm{diag}( \mathbf{\nu} )
\end{equation}
where $\mathbf{\nu}$ is an $N$-dimensional vector whose entries
are drawn from a uniform distribution between 0 and 1.

\subsection{Ensemble Kalman Filters (EnKFs)}
\label{enkf}

The EnKF is a method that replaces a single forecast state with an
ensemble of states. The spread of the ensemble about its mean gives an
approximation of the background error covariance and forecast
uncertainty, while the ensemble average gives the best guess of the
forecast. The EnKF was first introduced by \cite{evensen1994}.
For a comprehensive overview of ensemble filters, see
\cite{evensen2003}. It was shown that if the observation, which has
random error with covariance $\mathbf{R}$, is perturbed with $P$
random errors (again with covariance $\mathbf{R}$), to make an
$P$-member ensemble of independent observations $\left\{ \yvec_i
\right\}$, then the background error covariance can be written
\citep{evensen2003}
\begin{equation} \mathbf{B} \approx \frac{1}{P-1} \sum_{i=1}^{P}
  (\xvec_i^b - \overline{\xvec^b})(\xvec_i^b - \overline{\xvec^b})^T =
  \frac{1}{P-1} \mathbf{X}^b {\mathbf{X}^b}^T
\end{equation} 
which is simply the unbiased average outer product of background
perturbations $\mathbf{X}^b=[\xvec'^b_1,\ldots, \xvec'^b_P]$. The
background forecast of ensemble member $i$ is denoted $\xvec_i^b$,
$\overline{\xvec^b}$ is the background forecast ensemble average, and
$\xvec'^b_i = \xvec^b_i - \overline{\xvec^b}$ is the $i$th member's
deviation from the mean. In this case, each ensemble member is updated
according to the KF equations for their associated observation

The advantages of the EnKF are many: there is no linear tangent model
to compute, the number of ensemble members can be small
($\mathcal{O}(10^2)$ for NWP) relative to the dimensionality of the
state space, and prior knowledge about the structure of the forecast
errors is not necessary. Currently, 4D-Var (like 3D-Var but also
taking into account older observations) and ensemble filters are the
most promising candidates being considered to replace 3D-Var in
operational NWP.

As with the EKF, ensemble filters tend to underestimate the background
error, resulting in an ensemble spread which is typically less 
too small. We again used multiplicative inflation of the 
background error, a common method shown to be successful in
\cite{evensen2003,whitaker2002,annan2004,yang2006, kalnay2007}.
This is accomplished by setting
\begin{equation}
  \label{eq:infmultenkf}
  \mathbf{X}^b \leftarrow (1+\Delta)^{1/2}  \mathbf{X}^b 
\end{equation}
before the analysis step.
Additive inflation proved crucial to stabilizing both EnKFs tested.
Without it, the filters sometimes worked but only with $\Delta \gg 1$;
$\Delta$ is supposed to be a small parameter. As in the EKF,
additive inflation is applied immediately after the analysis
step, but in this case the noise is added to the analysis ensemble
states
\begin{equation}
  \label{eq:infaddenkf}
  \xvec_i^a \leftarrow \xvec_i^a + \mu \, \nu
\end{equation}
for all $i = 1, \ldots, P$. The noise $\nu$ is, again, an 
$N$-dimensional random vector with entries drawn from the
uniform distribution between 0 and 1.

\subsection{Ensemble Square Root Filter (EnSRF)} 
The original EnKF adds noise to create linearly
independent observations and is classified as a \textit{perturbed
  observations} method \citep{kalnay2007}. This necessarily introduces
additional sampling error into the forecast. For this reason,
\cite{whitaker2002} introduced the ensemble square root
filter (EnSRF) as an improved EnKF. In the EnSRF, the ensemble mean is
updated with the traditional Kalman gain (Eqn.~\eqref{eq:KFgain})
\begin{align} \overline{\xvec^a} &= \overline{\xvec^b}+\mathbf{K}
  (\yvec-\mathcal{H}\overline{\xvec^b})
\end{align} and deviations from the mean are updated by
\begin{align} \mathbf{X}^a &= (1-\tilde{\mathbf{K}}
  \mathcal{H})\mathbf{X}^b
\end{align} where
\begin{equation}
  \label{
    eq:EnSRFmult} 
  \tilde{\mathbf{K}} =
  \mathbf{BH}^T \left[ \left( \sqrt{\mathbf{HBH}^T+\mathbf{R}}
    \right)^{-1}\right]^T 
  \times
  \left[\sqrt{\mathbf{HBH}^T+\mathbf{R}}+\sqrt{\mathbf{R}} \right]^{-1}
  \;.
\end{equation}
When the observation is a scalar, it can be shown that
\begin{align}\label{eq:EnSRFscalar} \tilde{\mathbf{K}} &= \left( 1 +
    \sqrt{\frac{\mathbf{R}}{\mathbf{HBH}^T+\mathbf{R}}} \right)^{-1}
  \mathbf{K} \;.
\end{align} If observation errors are uncorrelated ($\mathbf{R}$ is
diagonal), then Eqn.~\eqref{eq:EnSRFscalar} can be used to process
observations one at a time \citep{whitaker2002}. The updated analysis
ensemble is then $\{\xvec_i^a\}$, where
$\xvec_i^a=\overline{\xvec^a}+\xvec_i'^a$. Square root filters have
better numerical stability and speed than their standard KF
counterparts. The Potter square root filter was
employed for navigation in the Lunar Module of the Apollo program
\citep{savely1972a}.

\subsection{Ensemble Transform Kalman Filter (ETKF)}

The ETKF is another type of deterministic square root filter. In this
variant, the analysis perturbations are assumed to be equal to the
background perturbations postmultiplied by a transformation matrix
$\mathbf{T}$ so that the analysis error covariance satisfies
Eqn.~\eqref{eq:KFanaylsiscov}. The analysis covariance is written
\[ \mathbf{A}=\frac{1}{P-1}\mathbf{X}^a {\mathbf{X}^a}^T= \mathbf{X}^b
\hat{\mathbf{A}} {\mathbf{X}^b}^T
\] where
$\hat{\mathbf{A}}=[(P-1)\mathbf{I}+(\mathbf{HX}^b)^T\mathbf{R}^{-1}(\mathbf{HX}^b)]^{-1}$
. The analysis perturbations are $\mathbf{X}^a=\mathbf{X}^b
\mathbf{T}$, where $\mathbf{T}=[(P-1)\hat{\mathbf{A}}]^{1/2}$. See
\cite{kalnay2007} for further details.

The local ensemble transform filter (LETKF) is a variant that computes
the analysis at a given gridpoint using only local observations. This
allows for efficient parallelization. Localization removes spurious
long-distance correlations from $\mathbf{B}$ and allows greater
flexibility in the global analysis by allowing different linear
combinations of ensemble members at different spatial locations
\citep{kalnay2007,hunt2007}.

\subsection{Tuning Parameters}
\begin{figure}
  \centering
  \includegraphics[width=\linewidth]{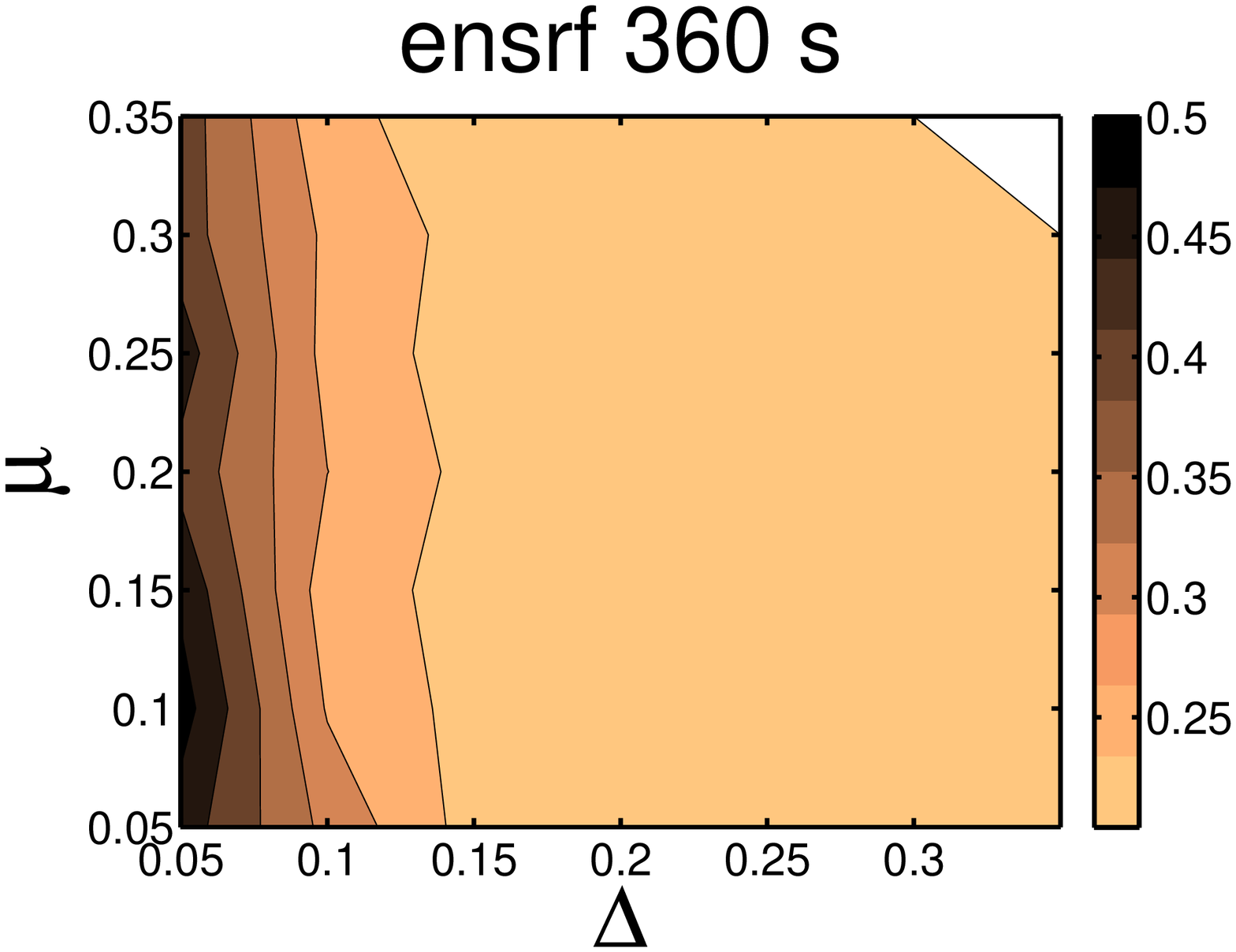}
  \caption{Contour plot of forecast error using the
  EnSRF and an assimilation window of 360 s, as the inflation
  parameters $\Delta$ and $\mu$ were varied. The chosen parameters
  $\Delta = 0.15$ and $\mu =0.25$ lie in the region of lowest error.}
  \label{fig:S2}
\end{figure}

Table~\ref{tbl:datuneparam} lists the tuning parameters used
for the DA experiments. The tuning was done manually.
Sensitivity of model error to the tuning parameters was checked
by creating a course contour plot of background error
for assimilation windows of 2, 4, 6, 8, and 10 minutes for 
each filter; an example is shown in Fig.~\ref{fig:S2}.

\begin{table}
  \begin{center}
    \begin{tabular}{l|l|l}
      \multicolumn{3}{|c|}{EnSRF \& ETKF} \\
      \hline
      analysis window (s) & $\Delta$ & $\mu$ \\
      \hline\hline
      30  & 0.15 & 0.25 \\
      60  & 0.15 & 0.25 \\
      90  & 0.15 & 0.25 \\
      120 & 0.15 & 0.25 \\
      150 & 0.15 & 0.25 \\
      180 & 0.15 & 0.25 \\
      210 & 0.15 & 0.25 \\
      240 & 0.15 & 0.25 \\
      270 & 0.15 & 0.25 \\
      300 & 0.15 & 0.25 \\
      330 & 0.15 & 0.25 \\
      360 & 0.15 & 0.25 \\
      390 & 0.2 & 0.25 \\
      420 & 0.2 & 0.25 \\
      450 & 0.2 & 0.25 \\
      480 & 0.25 & 0.25 \\
      510 & 0.25 & 0.25 \\
      540 & 0.25 & 0.25 \\
      570 & 0.25 & 0.25 \\
      600 & 0.25 & 0.25
    \end{tabular}
    \hspace{1em}
    \begin{tabular}{l|l|l}
      \multicolumn{3}{|c|}{EKF} \\
      \hline
      analysis window (s) & $\Delta$ & $\mu$ \\
      \hline\hline
      30  & 0.15 & 0.25  \\
      60  & 0.15 & 0.25  \\
      90  & 0.15 & 0.25  \\
      120 & 0.15 & 0.25  \\
      150 & 0.15 & 0.25  \\
      180 & 0.15 & 0.25  \\
      210 & 0.15 & 0.25  \\
      240 & 0.15 & 0.25  \\
      270 & 0.15 & 0.25  \\
      300 & 0.15 & 0.25  \\
      330 & 0.15 & 0.25  \\
      360 & 0.15 & 0.25 \\
      390 & 0.2 & 0.25  \\
      420 & 0.2 & 0.25  \\
      450 & 0.2 & 0.25  \\
      480 & 0.2 & 0.25  \\
      510 & 0.2 & 0.25  \\
      540 & 0.25 & 0.25  \\
      570 & 0.25 & 0.25  \\
      600 & 0.25 & 0.25 
    \end{tabular}
  \caption{Inflation tuning parameters used in DA experiments.}
  \label{tbl:datuneparam}
\end{center}
\end{table}

\section{Flow reversal forecast tuning}
\begin{figure}
  \centering
  \includegraphics[width=\linewidth]{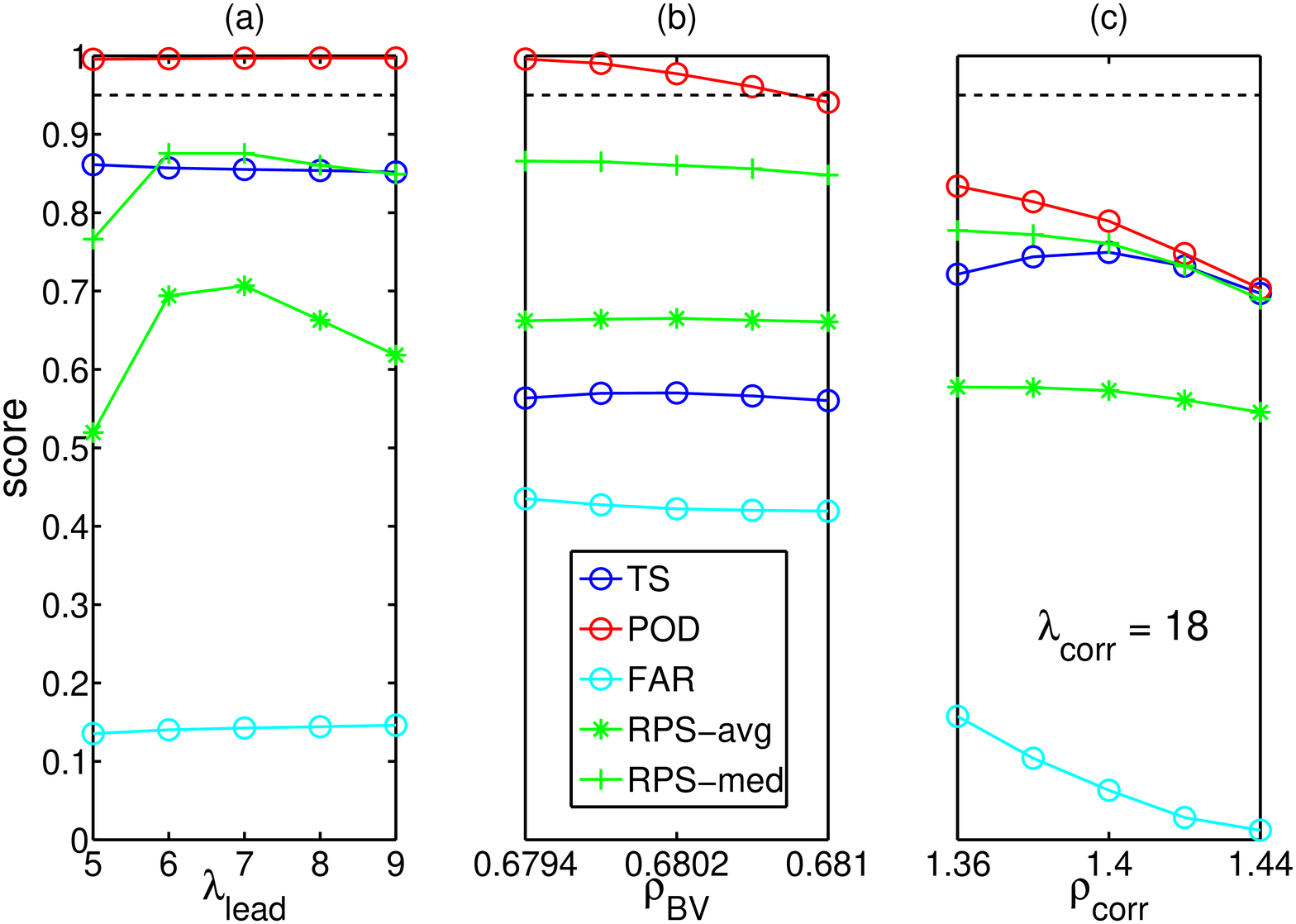}
  
  \caption{Tuning the flow reversal forecasts methods by varying
    (a) $\rho_\mathrm{corr}$ for the correlation test, 
    (b) $\rho_\mathrm{BV}$ for the BV test, and 
    (c) $\lambda_\mathrm{lead}$ for lead test. The dashed
    line is at 95\%.
    There are tradeoffs among the various skill scores 
    (POD, TS, RPS-avg, RPS-med) and costs (FAR). 
    The chosen, final parameters 
    $\lambda_\text{lead} = 7$, 
    $\rho_\text{BV} = 0.6802$, and 
    $\rho_\text{corr} = 1.4$ (with $\lambda_\text{corr} = 18$)
    appear in the middle of each plot.
    \label{fig:regchgtune}
  }
\end{figure}

In Fig.~\ref{fig:regchgtune},
we present skill score curves as the
tuning parameters of the three flow reversal forecasts are varied.
We used a sequence of plots of this type to inform our
tuning of the various flow reversal tests.

\begin{figure}
  \centering

  
  \includegraphics[width=.7\linewidth, totalheight=.75\textheight,
  keepaspectratio=true]{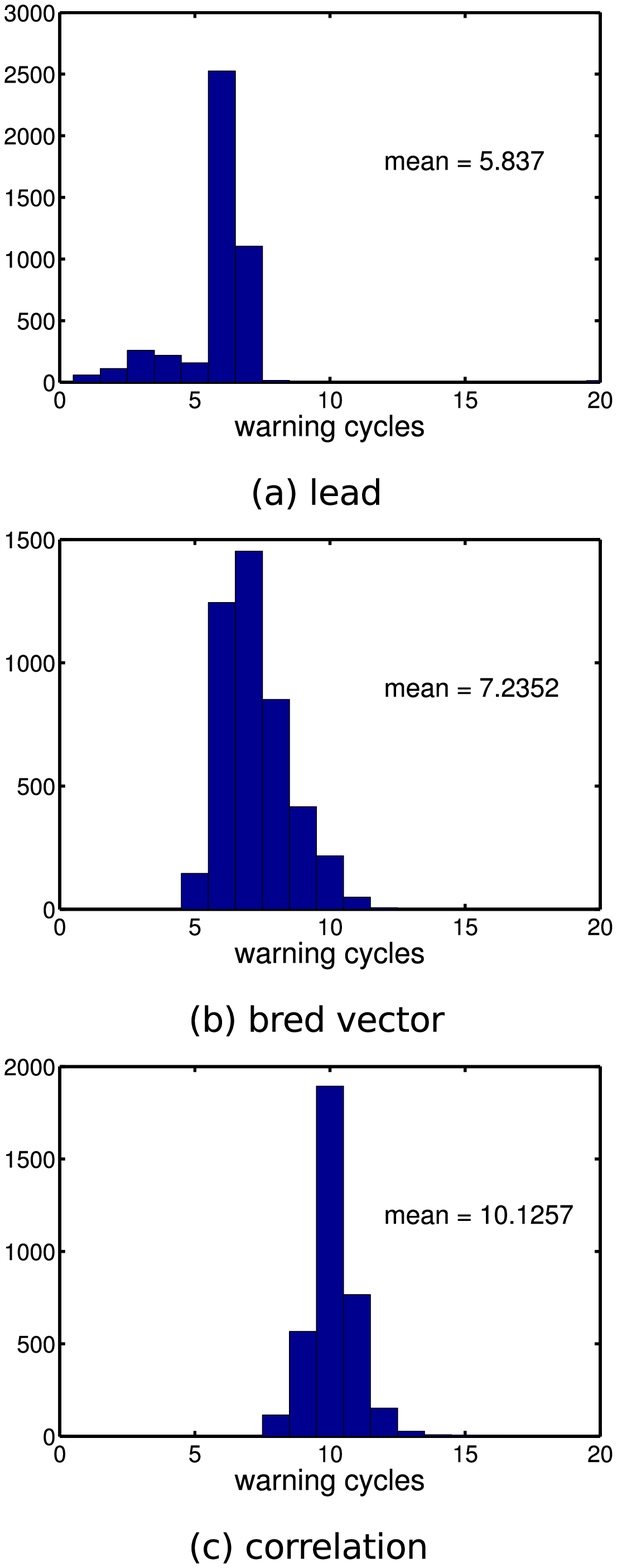}

  \caption{Warning time distributions for the three flow reversal tests,
    in units of 30 s analysis cycles. The warning time is the time interval
    between the test being triggered and observation of the the actual
    flow reversal. The correlation test gives the
    earliest average warning time, followed by the BV and lead tests.}
  \label{fig:warntime}
\end{figure}

In Fig.~\ref{fig:warntime}, we also present warning time 
histograms for the different tests. If early detection of flow
reversals is desirable, then the warning times
tell us how the tests compare.

\label{sec:revoccur.skill}

\end{document}